\documentstyle[11pt]{article}
\title{The Intrinsic Normal Cone}
\author{K. Behrend and B. Fantechi}

\catcode`\@=11

\font\tenmsx=msam10
\font\sevenmsx=msam7
\font\fivemsx=msam5
\font\tenmsy=msbm10
\font\sevenmsy=msbm7
\font\fivemsy=msbm5
\newfam\msxfam
\newfam\msyfam
\textfont\msxfam=\tenmsx  \scriptfont\msxfam=\sevenmsx
  \scriptscriptfont\msxfam=\fivemsx
\textfont\msyfam=\tenmsy  \scriptfont\msyfam=\sevenmsy
  \scriptscriptfont\msyfam=\fivemsy

\def\hexnumber@#1{\ifcase#1 0\or1\or2\or3\or4\or5\or6\or7\or8\or9\or
        A\or B\or C\or D\or E\or F\fi }

\font\teneuf=eufm10
\font\seveneuf=eufm7
\font\fiveeuf=eufm5
\newfam\euffam
\textfont\euffam=\teneuf
\scriptfont\euffam=\seveneuf
\scriptscriptfont\euffam=\fiveeuf
\def\frak{\ifmmode\let\next\frak@\else
 \def\next{\Err@{Use \string\frak\space only in math mode}}\fi\next}
\def\goth{\relaxnext@\ifmmode\let\next\frak@\else
 \def\next{\Err@{Use \string\goth\space only in math mode}}\fi\next}
\def\frak@#1{{\frak@@{#1}}}
\def\frak@@#1{\fam\euffam#1}
\edef\msx@{\hexnumber@\msxfam}
\edef\msy@{\hexnumber@\msyfam}

\mathchardef\boxdot="2\msx@00
\mathchardef\boxplus="2\msx@01
\mathchardef\boxtimes="2\msx@02
\mathchardef\square="0\msx@03
\mathchardef\blacksquare="0\msx@04
\mathchardef\centerdot="2\msx@05
\mathchardef\lozenge="0\msx@06
\mathchardef\blacklozenge="0\msx@07
\mathchardef\circlearrowright="3\msx@08
\mathchardef\circlearrowleft="3\msx@09
\mathchardef\rightleftharpoons="3\msx@0A
\mathchardef\leftrightharpoons="3\msx@0B
\mathchardef\boxminus="2\msx@0C
\mathchardef\Vdash="3\msx@0D
\mathchardef\Vvdash="3\msx@0E
\mathchardef\vDash="3\msx@0F
\mathchardef\twoheadrightarrow="3\msx@10
\mathchardef\twoheadleftarrow="3\msx@11
\mathchardef\leftleftarrows="3\msx@12
\mathchardef\rightrightarrows="3\msx@13
\mathchardef\upuparrows="3\msx@14
\mathchardef\downdownarrows="3\msx@15
\mathchardef\upharpoonright="3\msx@16

\mathchardef\downharpoonright="3\msx@17
\mathchardef\upharpoonleft="3\msx@18
\mathchardef\downharpoonleft="3\msx@19
\mathchardef\rightarrowtail="3\msx@1A
\mathchardef\leftarrowtail="3\msx@1B
\mathchardef\leftrightarrows="3\msx@1C
\mathchardef\rightleftarrows="3\msx@1D
\mathchardef\Lsh="3\msx@1E
\mathchardef\Rsh="3\msx@1F
\mathchardef\rightsquigarrow="3\msx@20
\mathchardef\leftrightsquigarrow="3\msx@21
\mathchardef\looparrowleft="3\msx@22
\mathchardef\looparrowright="3\msx@23
\mathchardef\circeq="3\msx@24
\mathchardef\succsim="3\msx@25
\mathchardef\gtrsim="3\msx@26
\mathchardef\gtrapprox="3\msx@27
\mathchardef\multimap="3\msx@28
\mathchardef\therefore="3\msx@29
\mathchardef\because="3\msx@2A
\mathchardef\doteqdot="3\msx@2B

\mathchardef\triangleq="3\msx@2C
\mathchardef\precsim="3\msx@2D
\mathchardef\lesssim="3\msx@2E
\mathchardef\lessapprox="3\msx@2F
\mathchardef\eqslantless="3\msx@30
\mathchardef\eqslantgtr="3\msx@31
\mathchardef\curlyeqprec="3\msx@32
\mathchardef\curlyeqsucc="3\msx@33
\mathchardef\preccurlyeq="3\msx@34
\mathchardef\leqq="3\msx@35
\mathchardef\leqslant="3\msx@36
\mathchardef\lessgtr="3\msx@37
\mathchardef\backprime="0\msx@38
\mathchardef\risingdotseq="3\msx@3A
\mathchardef\fallingdotseq="3\msx@3B
\mathchardef\succcurlyeq="3\msx@3C
\mathchardef\geqq="3\msx@3D
\mathchardef\geqslant="3\msx@3E
\mathchardef\gtrless="3\msx@3F
\mathchardef\sqsubset="3\msx@40
\mathchardef\sqsupset="3\msx@41
\mathchardef\vartriangleright="3\msx@42
\mathchardef\vartriangleleft="3\msx@43
\mathchardef\trianglerighteq="3\msx@44
\mathchardef\trianglelefteq="3\msx@45
\mathchardef\bigstar="0\msx@46
\mathchardef\between="3\msx@47
\mathchardef\blacktriangledown="0\msx@48
\mathchardef\blacktriangleright="3\msx@49
\mathchardef\blacktriangleleft="3\msx@4A
\mathchardef\vartriangle="0\msx@4D
\mathchardef\blacktriangle="0\msx@4E
\mathchardef\triangledown="0\msx@4F
\mathchardef\eqcirc="3\msx@50
\mathchardef\lesseqgtr="3\msx@51
\mathchardef\gtreqless="3\msx@52
\mathchardef\lesseqqgtr="3\msx@53
\mathchardef\gtreqqless="3\msx@54
\mathchardef\Rrightarrow="3\msx@56
\mathchardef\Lleftarrow="3\msx@57
\mathchardef\veebar="2\msx@59
\mathchardef\barwedge="2\msx@5A
\mathchardef\doublebarwedge="2\msx@5B
\mathchardef\angle="0\msx@5C
\mathchardef\measuredangle="0\msx@5D
\mathchardef\sphericalangle="0\msx@5E
\mathchardef\varpropto="3\msx@5F
\mathchardef\smallsmile="3\msx@60
\mathchardef\smallfrown="3\msx@61
\mathchardef\Subset="3\msx@62
\mathchardef\Supset="3\msx@63
\mathchardef\Cup="2\msx@64

\mathchardef\Cap="2\msx@65

\mathchardef\curlywedge="2\msx@66
\mathchardef\curlyvee="2\msx@67
\mathchardef\leftthreetimes="2\msx@68
\mathchardef\rightthreetimes="2\msx@69
\mathchardef\subseteqq="3\msx@6A
\mathchardef\supseteqq="3\msx@6B
\mathchardef\bumpeq="3\msx@6C
\mathchardef\Bumpeq="3\msx@6D
\mathchardef\lll="3\msx@6E

\mathchardef\ggg="3\msx@6F

\mathchardef\circledS="0\msx@73
\mathchardef\pitchfork="3\msx@74
\mathchardef\dotplus="2\msx@75
\mathchardef\backsim="3\msx@76
\mathchardef\backsimeq="3\msx@77
\mathchardef\complement="0\msx@7B
\mathchardef\intercal="2\msx@7C
\mathchardef\circledcirc="2\msx@7D
\mathchardef\circledast="2\msx@7E
\mathchardef\circleddash="2\msx@7F
\def\ulcorner{\delimiter"4\msx@70\msx@70 }
\def\urcorner{\delimiter"5\msx@71\msx@71 }
\def\llcorner{\delimiter"4\msx@78\msx@78 }
\def\lrcorner{\delimiter"5\msx@79\msx@79 }
\def\yen{\mathhexbox\msx@55 }
\def\checkmark{\mathhexbox\msx@58 }
\def\circledR{\mathhexbox\msx@72 }
\def\maltese{\mathhexbox\msx@7A }
\mathchardef\lvertneqq="3\msy@00
\mathchardef\gvertneqq="3\msy@01
\mathchardef\nleq="3\msy@02
\mathchardef\ngeq="3\msy@03
\mathchardef\nless="3\msy@04
\mathchardef\ngtr="3\msy@05
\mathchardef\nprec="3\msy@06
\mathchardef\nsucc="3\msy@07
\mathchardef\lneqq="3\msy@08
\mathchardef\gneqq="3\msy@09
\mathchardef\nleqslant="3\msy@0A
\mathchardef\ngeqslant="3\msy@0B
\mathchardef\lneq="3\msy@0C
\mathchardef\gneq="3\msy@0D
\mathchardef\npreceq="3\msy@0E
\mathchardef\nsucceq="3\msy@0F
\mathchardef\precnsim="3\msy@10
\mathchardef\succnsim="3\msy@11
\mathchardef\lnsim="3\msy@12
\mathchardef\gnsim="3\msy@13
\mathchardef\nleqq="3\msy@14
\mathchardef\ngeqq="3\msy@15
\mathchardef\precneqq="3\msy@16
\mathchardef\succneqq="3\msy@17
\mathchardef\precnapprox="3\msy@18
\mathchardef\succnapprox="3\msy@19
\mathchardef\lnapprox="3\msy@1A
\mathchardef\gnapprox="3\msy@1B
\mathchardef\nsim="3\msy@1C
\mathchardef\ncong="3\msy@1D

\mathchardef\varsubsetneq="3\msy@20
\mathchardef\varsupsetneq="3\msy@21
\mathchardef\nsubseteqq="3\msy@22
\mathchardef\nsupseteqq="3\msy@23
\mathchardef\subsetneqq="3\msy@24
\mathchardef\supsetneqq="3\msy@25
\mathchardef\varsubsetneqq="3\msy@26
\mathchardef\varsupsetneqq="3\msy@27
\mathchardef\subsetneq="3\msy@28
\mathchardef\supsetneq="3\msy@29
\mathchardef\nsubseteq="3\msy@2A
\mathchardef\nsupseteq="3\msy@2B
\mathchardef\nparallel="3\msy@2C
\mathchardef\nmid="3\msy@2D
\mathchardef\nshortmid="3\msy@2E
\mathchardef\nshortparallel="3\msy@2F
\mathchardef\nvdash="3\msy@30
\mathchardef\nVdash="3\msy@31
\mathchardef\nvDash="3\msy@32
\mathchardef\nVDash="3\msy@33
\mathchardef\ntrianglerighteq="3\msy@34
\mathchardef\ntrianglelefteq="3\msy@35
\mathchardef\ntriangleleft="3\msy@36
\mathchardef\ntriangleright="3\msy@37
\mathchardef\nleftarrow="3\msy@38
\mathchardef\nrightarrow="3\msy@39
\mathchardef\nLeftarrow="3\msy@3A
\mathchardef\nRightarrow="3\msy@3B
\mathchardef\nLeftrightarrow="3\msy@3C
\mathchardef\nleftrightarrow="3\msy@3D
\mathchardef\divideontimes="2\msy@3E
\mathchardef\varnothing="0\msy@3F
\mathchardef\nexists="0\msy@40
\mathchardef\mho="0\msy@66
\mathchardef\eth="0\msy@67
\mathchardef\eqsim="3\msy@68
\mathchardef\beth="0\msy@69
\mathchardef\gimel="0\msy@6A
\mathchardef\daleth="0\msy@6B
\mathchardef\lessdot="3\msy@6C
\mathchardef\gtrdot="3\msy@6D
\mathchardef\ltimes="2\msy@6E
\mathchardef\rtimes="2\msy@6F
\mathchardef\shortmid="3\msy@70
\mathchardef\shortparallel="3\msy@71
\mathchardef\smallsetminus="2\msy@72
\mathchardef\thicksim="3\msy@73
\mathchardef\thickapprox="3\msy@74
\mathchardef\approxeq="3\msy@75
\mathchardef\succapprox="3\msy@76
\mathchardef\precapprox="3\msy@77
\mathchardef\curvearrowleft="3\msy@78
\mathchardef\curvearrowright="3\msy@79
\mathchardef\digamma="0\msy@7A
\mathchardef\varkappa="0\msy@7B
\mathchardef\hslash="0\msy@7D
\mathchardef\hbar="0\msy@7E
\mathchardef\backepsilon="3\msy@7F
\def\Bbb{\ifmmode\let\next\Bbb@\else
 \def\next{\errmessage{Use \string\Bbb\space only in math mode}}\fi\next}
\def\Bbb@#1{{\Bbb@@{#1}}}
\def\Bbb@@#1{\fam\msyfam#1}

\catcode`\@=12

\newtheorem{prop}{Proposition}[section]
\newtheorem{lem}[prop]{Lemma}

\newtheorem{cor}[prop]{Corollary}
\newtheorem{them}[prop]{Theorem}

\newtheorem{defnp}[prop]{Definition}
\newtheorem{numconp}[prop]{Construction}
\newtheorem{numexp}[prop]{Example}
\newtheorem{numrmkp}[prop]{Remark}
\newtheorem{numrmksp}[prop]{Remarks}
\newtheorem{conp}[prop]{}

\newenvironment{defn}{\begin{defnp}\rm}{\end{defnp}}

\newenvironment{numex}{\begin{numexp}\rm}{\end{numexp}}
\newenvironment{numrmk}{\begin{numrmkp}\rm}{\end{numrmkp}}

\newtheorem{warningp}{Warning}
\newtheorem{notep}{Note}
\newtheorem{claimp}{Claim}
\newtheorem{examplep}{Example}
\newtheorem{examplesp}{Examples}
\newtheorem{rmkp}{Remark}
\newtheorem{rmksp}{Remarks}

\newenvironment{note}{\begin{notep}\rm}{\end{notep}}

\newenvironment{example}{\begin{examplep}\rm}{\end{examplep}}
\newenvironment{examples}{\begin{examplesp}\rm}{\end{examplesp}}
\newenvironment{rmk}{\begin{rmkp}\rm}{\end{rmkp}}

\newenvironment{pf}{\begin{trivlist}\item[]{\sc Proof.}}%
            {\nolinebreak $\Box$ \end{trivlist}}

\newcommand{\noprint}[1]{}

\renewcommand{\tilde}{\widetilde}
\newcommand{\dual}[1]{{{#1}^{\vee}}}

\newcommand{\qed}{{\nolinebreak $\Box$}}

\newcommand{\et}{{\mbox{\tiny \'{e}t}}}

\newcommand{\fl}{{\mbox{\tiny fl}}}

\newcommand{\upst}{^{\ast}}

\newcommand{\upsh}{^{!}}
\newcommand{\lst}{_{\ast}}

\newcommand{\com}{^{\scriptscriptstyle\bullet}}

\newcommand{\argument}{{{\,\cdot\,}}}

\newcommand{\DD}{{\frak D}}

\newcommand{\EE}{{\frak E}}
\newcommand{\FF}{{\frak F}}

\newcommand{\NN}{{\frak N}}

\newcommand{\Mm}{{\frak m}}

\newcommand{\CC}{{\frak C}}
\newcommand{\zz}{{\Bbb Z}}

\newcommand{\aaa}{{\Bbb A}}

\newcommand{\qq}{{\Bbb Q}}

\newcommand{\tT}{{\cal T}}

\renewcommand{\O}{{\cal O}}

\newcommand{\eE}{{\cal E}}
\newcommand{\iI}{{\cal I}}
\newcommand{\jJ}{{\cal J}}
\newcommand{\gG}{{\cal G}}
\newcommand{\fF}{{\cal F}}

\newcommand{\hH}{{\cal H}}
\newcommand{\sS}{{\cal S}}
\newcommand{\pP}{{\cal P}}

\newcommand{\del}{\partial}
\newcommand{\resto}{{ \mid }}

\newcommand{\rk}{\mathop{\rm rk}}
\newcommand{\tot}{\mathop{\rm tot}}

\newcommand{\td}{\mathop{\rm td}}

\newcommand{\pr}{\mathop{\rm pr}\nolimits}

\newcommand{\Mod}{\mathop{\rm Mod}\nolimits}
\newcommand{\Mor}{\mathop{\rm Mor}\nolimits}

\newcommand{\im}{\mathop{\rm im}}
\newcommand{\ob}{\mathop{\rm ob}}

\newcommand{\cok}{\mathop{\rm cok}}

\newcommand{\spec}{\mathop{\rm Spec}\nolimits}
\newcommand{\Spec}{\spec}
\newcommand{\Sym}{\mathop{\rm Sym}\nolimits}

\newcommand{\id}{\mathop{\rm id}}
\newcommand{\Hom}{\mathop{\rm Hom}\nolimits}
\newcommand{\Ext}{\mathop{\rm Ext}\nolimits}

\newcommand{\ltensor}{\mathop{\displaystyle\stackrel{L}{\otimes}}}
\newcommand{\sheafhom}{\mathop{\rm {\mit{ \hH\! om}}}\nolimits}

\newcommand{\sheaftor}{\mathop{\rm {\mit{ \tT\! or}}}\nolimits}

\newcommand{\Aut}{\mathop{\rm Aut}\nolimits}

\newcommand{\comp}{\mathbin{{\scriptstyle\circ}}}
\newcommand{\ol}{\overline}
\newcommand{\ul}{\underline}

\newcommand{\ldiag}[1]%
       {\makebox[0cm]{${\scriptstyle#1}\downarrow\phantom{\scriptstyle#1}$}}
\newcommand{\ldiagup}[1]%
       {\makebox[0cm]{${\scriptstyle#1}\uparrow\phantom{\scriptstyle#1}$}}
\newcommand{\rdiag}[1]%
       {\makebox[0cm]{$\phantom{\scriptstyle#1}\downarrow{\scriptstyle#1}$}}
\newcommand{\sediagr}[1]%
       {\makebox[0cm]{$\phantom{\scriptstyle#1}\searrow{\scriptstyle#1}$}}
\newcommand{\nediagr}[1]%
       {\makebox[0cm]{$\phantom{\scriptstyle#1}\nearrow{\scriptstyle#1}$}}
\newcommand{\rdiagup}[1]%
       {\makebox[0cm]{$\phantom{\scriptstyle#1}\uparrow{\scriptstyle#1}$}}
\newcommand{\swdiag}[1]%
       {\makebox[0cm]{$\phantom{\scriptstyle#1}\swarrow{\scriptstyle#1}$}}
\newcommand{\sediag}[1]%
       {\makebox[0cm]{${\scriptstyle#1}\searrow\phantom{\scriptstyle#1}$}}
\newcommand{\nediag}[1]%
       {\makebox[0cm]{${\scriptstyle#1}\nearrow\phantom{\scriptstyle#1}$}}

\newcommand{\longiso}{\stackrel{\textstyle\sim}{\longrightarrow}}
\newcommand{\iso}{\stackrel{\sim}{\rightarrow}}

\newcommand{\comdia}[9]{%
\begin{array}{ccc}
#1 & \stackrel{#2}{\longrightarrow} & #3 \\
\ldiag{#4} & #5 & \rdiag{#6} \\
#7 & \stackrel{#8}{\longrightarrow} & #9
\end{array}}

\newcommand{\comtri}[6]{%
\begin{array}{ccc}
#1 & \stackrel{#2}{\longrightarrow} & #3         \\
   & \sediag{#4}                    & \rdiag{#5} \\
   &                                & #6
\end{array}}

\newcommand{\inversecomtri}[6]{%
\begin{array}{ccc}
#1           &                                &            \\
\ldiag{#2}   & \sediagr{#3}                   &            \\
#4           & \stackrel{#5}{\longrightarrow} & #6
\end{array}}

\newcommand{\overtoparrow}%
{\makebox[0cm]{\beginpicture
\setcoordinatesystem units <.8cm,.4cm> point at 0 0
\setplotarea x from -3 to 3, y from 0 to 1
\setquadratic
\plot -3 0 0 1 3 0 /
\put{\vector(3,-1){0}}[Bl] at 3 0
\endpicture}}

\newcommand{\underbottomarrow}%
{\makebox[0cm]{\beginpicture
\setcoordinatesystem units <.8cm,.4cm> point at 0 0
\setplotarea x from -3 to 3, y from 0 to 1
\setquadratic
\plot -3 1 0 0 3 1 /
\put{\vector(3,1){0}}[Bl] at 3 1
\endpicture}}

\newcommand{\ses}[5]%
{0\longrightarrow#1\stackrel{#2}{ \longrightarrow}#3\stackrel{#4}{
\longrightarrow}#5\longrightarrow0}

\newcommand{\dt}[6]%
{#1\stackrel{#2}{\longrightarrow}#3 \stackrel{#4}{\longrightarrow}#5
\stackrel{#6}{\longrightarrow} #1[1]}

\begin{document} \sloppy
\date{January 13, 1996}
\maketitle
\begin{abstract}
We suggest a  construction of  virtual fundamental classes
of certain types of moduli spaces.
\end{abstract}
\tableofcontents
\section{Introduction}
Moduli spaces in algebraic geometry often have an expected dimension
at each point, which is a lower bound for the dimension at that
point. For instance, the moduli space of smooth, complex projective
$n$-dimensional varieties with ample canonical class has expected
dimension $h^1(V,T_V)-h^2(V,T_V)$ at a point $[V]$. In general, the
expected dimension will vary with the point; however, in some
significant cases it will stay constant on connected components. In
the previous example, this is the case if $n\le 2$, for then the
expected dimension is $-\chi(V,T_V)$. In some cases the dimension
coincides with the expected dimension, in others it does so under some
genericity assumptions. However, it can happen that there is no way to
get a space of the expected dimension; it is also possible that
special cases with bigger dimension are easier to understand and to
deal with than the generic case.

When we have a moduli space $X$ which has a well-defined expected
dimension, it can be useful to be able to construct in its Chow ring a
class of the expected dimension. The main examples we have in mind are
Donaldson theory (with $X$ the moduli space of torsion-free,
semi-stable sheaves on a surface) and the Gromov-Witten invariants
(with $X$ the moduli space of stable maps from curves of genus $g$ to
a fixed projective variety).  In this paper we deal with the problem
of defining such a class in a very general set-up; the construction is
divided into two steps.

First, given any Deligne-Mumford stack $X$, we associate to it an
algebraic stack $\CC_X$ over $X$ of pure dimension zero, its {\em
intrinsic normal cone}. This has nothing to do with $X$ being a moduli
space; it is just an intrinsic invariant, whose structure is related
to the singularities of $X$ (see for instance Proposition
\ref{colci}).

Then, we define the concept of an obstruction theory and of a perfect
obstruction theory for $X$. To say that $X$ has an obstruction theory
means, very roughly speaking, that we are given locally on $X$ an
(equivalence class of) morphisms of vector bundles such that at each
point the kernel of the induced linear map of vector spaces is the
tangent space to $X$, and the cokernel is a space of
obstructions. Usually, if $X$ is a moduli space then it has an
obstruction theory, and if this is perfect then the expected dimension
is constant on $X$.  Once we are given an obstruction theory, with the
additional (technical) assumption that it admits a global resolution,
we can define a virtual fundamental class of the expected dimension.

An application of the results of this work is contained in a paper
\cite{gwi} by the first author. There Gromov-Witten invariants are
constructed for any genus, any target variety and the axioms listed in
\cite{BM} are verified.

We now give a more detailed outline of the contents of the paper. In
the first section we recall what we need about cones and we introduce
the notion of cone stacks over a Deligne-Mumford stack $X$. These are
Artin stacks which are locally the quotient of a cone by a vector
bundle acting on it. We call a cone {\em abelian }if it is defined as
$\spec\Sym \fF$, where $\fF$ is a coherent sheaf on $X$. Every cone is
contained as a closed subcone in a minimal abelian one, which we call
its {\em abelian hull}. The notions of being abelian and of abelian hull
generalize immediately to cone stacks.

In the second section we
construct, for a complex $E\com$ in the derived category $D(\O_X)$
which satisfies some suitable assumptions (which we call Condition
($\star$), see Definition \ref{dost}), an associated abelian cone
stack $h^1/h^0((E\com\dual))$. In particular the cotangent complex
$L_X\com$ of $X$ satisfies Condition ($\star$), so we can define the
abelian cone stack $\NN_X:=h^1/h^0((L_X\com\dual))$, the {\em
intrinsic normal sheaf}.

The name is motivated in the third section, where $\NN_X$
is constructed more directly as follows: \'etale locally on $X$, embed
an open set $U$ of $X$ in a smooth scheme $W$, and take the stack
quotient of the normal sheaf (viewed as abelian cone) $N_{U/W}$ by the
natural action of $T_W|_U$. One can glue these abelian cone stacks
together to get $\NN_X$.  The intrinsic normal cone $\CC_X$ is the closed
subcone stack of $\NN_X$ defined by replacing
$N_{U/W}$ by the normal cone $C_{U/W}$ in the previous construction.

In the fourth section we describe the relationship between the intrinsic normal
sheaf of a Deligne-Mumford stack $X$ and the deformations of affine
$X$-schemes, showing in particular that $\NN_X$ carries obstructions for such
deformations. With this motivation, we introduce the notion of
obstruction theory for $X$. This is an object $E\com$ in the
derived category together with a morphism $E\com\to L_X\com$,
satisfying Condition ($\star$) and such that the induced
map $\NN_X\to h^1/h^0((E\com\dual))$ is a closed immersion.

An obstruction theory $E\com$ is called perfect if
$\EE=h^1/h^0((E\com\dual))$ is smooth over $X$. So we have a vector
bundle stack $\EE$ with a closed subcone stack $\CC_X$, and to define
the virtual fundamental class of $X$ with respect to $E\com$ we simply
intersect $\CC_X$ with the zero section of $\EE$. This construction
requires Chow groups for Artin stacks, which we do not have at our
disposal. There are several ways around  this problem. We choose to
assume that $E\com$ is globally given by a homomorphism of vector
bundles $F^{-1}\to F^0$. Then $\CC_X$ gives rise to a cone $C$ in
$F_1=\dual{F^{-1}}$ and we intersect $C$ with the zero section of
$F_1$.

Another approach, suggested by Kontsevich \cite{K}, is via virtual
structure sheaves (see Remark~\ref{vss}). The drawback of that
approach is that it requires a Riemann-Roch theorem for
Deligne-Mumford stacks, for which we do not know a reference.

In the sixth section we give some examples of how this construction
can be applied in some standard moduli problems. We consider the
following cases: a fiber of a morphism between smooth algebraic
stacks, the scheme of morphisms between two given projective schemes,
a moduli space for Gorenstein projective varieties.

In the seventh section we give a relative version of the intrinsic
normal cone and sheaf $\CC_{X/Y}$ and $\NN_{X/Y}$ for a morphism $X\to
Y$ with unramified diagonal of algebraic stacks; we are mostly
interested in the case where $Y$ is smooth and pure-dimensional, which
preserves many good properties of the absolute case (e.g., $\CC_{X/Y}$
is pure-dimensional). This is not needed in this paper, but will be
applied by the first author to give an algebraic definition of
Gromov-Witten classes for smooth projective varieties.

The starting point for this work was a talk by Jun Li at the AMS
Summer Institute on Algebraic Geometry, Santa Cruz 1995, where he
reported on joint work in progress with G.~Tian. Their construction,
in the complex analytic context, is based on the existence of the
Kuranishi map; by using it they define, under suitable assumptions, a
pure-dimensional cone in some bundle and get classes of the expected
dimension by intersecting with the zero section.

Our construction owes its existence to theirs; we started by trying to
understand and reformulate their results in an algebraic way, and
found stacks to be a convenient, intrinsic language. In our opinion
the introduction of stacks is very natural, and it seems almost
surprising that the intrinsic normal cone was not defined before. We
find it important to separate the construction of the cone, which can
be carried out for any Deligne-Mumford stack, from its embedding in a
vector bundle stack. We work completely in an algebraic context; of
course the whole paper could be rewritten without changes over the
category of analytic spaces.

\smallskip\noindent
{\em Acknowledgments. }This work was started in the inspiring
atmosphere of the Santa Cruz conference. A significant part of it was done
during
the authors' stay at the Max-Planck-Institut f\"ur Mathematik in Bonn, to which
both authors are grateful for hospitality and support. The second author is a
member of GNSAGA of CNR.

\subsection{Notations and Conventions}

Unless otherwise mentioned, we work over a fixed ground field $k$.

An {\em algebraic stack }is an algebraic stack over $k$ in the sense
of \cite{vdas} or \cite{laumon}. Unless mentioned otherwise, we assume
all algebraic stacks (in particular all algebraic spaces and all
schemes) to be quasi-separated and locally of finite type over $k$.

A {\em Deligne-Mumford stack }is an algebraic stack in the sense of
\cite{DM}, in other words an algebraic stack with unramified
diagonal. For a Deligne-Mumford stack $X$ we denote by $X_\fl$ the big
fppf-site and by $X_\et$ the small \'etale site of $X$. The associated
topoi of sheaves are denoted by the same symbols.

Recall that a complex of sheaves of modules is {\em of perfect
amplitude contained in $[a,b]$}, where $a,b\in\zz$, if, locally, it is
isomorphic (in the derived category) to a complex $E^a\to\ldots\to E^b$
of locally free sheaves of finite rank.

\section{Cones and Cone Stacks}

\subsection{Cones}

To fix notation we recall some basic facts about cones.

Let $X$ be a Deligne-Mumford stack. Let
\[S=\bigoplus_{i\geq0}S^i\]
be a graded quasi-coherent sheaf of $\O_X$-algebras such that
$S^0=\O_X$, $S^1$ is coherent and $S$ is generated locally by
$S^1$. Then the affine $X$-scheme $C=\spec S$ is called a {\em cone
}over $X$. A {\em morphism }of cones over $X$ is an $X$-morphism
induced by a graded morphism of graded sheaves of $\O_X$-algebras.  A
{\em closed subcone }is the image of a closed immersion of cones.  If
\[\begin{array}{ccc}
 & & C_2 \\
 & & \downarrow \\ C_1 & \longrightarrow & C_3
\end{array}\] is a diagram of cones over $X$, the fibered product
$C_1\times_{C_3}C_2$ is a cone over $X$.

Every cone $C\to X$ has a section $0:X\to C$, called the {\em vertex
}of $C$, and an $\aaa^1$-action (or a multiplicative contraction onto
the vertex), that is a morphism
\[\gamma:\aaa^1\times C\longrightarrow C\]
such that
\begin{enumerate}
\item \[\comtri{C}{(1,\id)}{\aaa^1\times C}{\id}{\gamma}{C}\] commutes,
\item \[\comtri{C}{(0,\id)}{\aaa^1\times C}{0}{\gamma}{C}\] commutes,
\item \[\comdia{\aaa^1\times\aaa^1\times C} {\id\times\gamma}
{\aaa^1\times C} {m\times\id}{}{\gamma} {\aaa^1\times C}{\gamma}{C}\]
commutes, where $m:\aaa^1\times\aaa^1\to\aaa^1$ is multiplication,
$m(x,y)=xy$.
\end{enumerate}
The vertex of $C$ is induced by the augmentation $S\to S^0$, the
$\aaa^1$-action is given by the grading of $S$. In fact, the morphism
$S\to S[x]$ giving rise to $\gamma$ maps $s\in S^i$ to $sx^i$.

Note that a morphism of cones is just a morphism respecting $0$ and
$\gamma$.

\subsection{Abelian Cones}

If $\fF$ is a coherent $\O_X$-module we get an associated cone
\[C(\fF)=\spec\Sym(\fF).\]
For any $X$-scheme $T$ we have
\[C(\fF)(T)=\Hom(\fF_T,\O_T),\]
so $C(\fF)$ is a group scheme over $X$. We call a cone of this form an
{\em abelian cone}. A fibered product of abelian cones is an abelian
cone. If $E$ is a vector bundle over $X$, then $E=C(\eE^{\vee})$,
where $\eE$ is the coherent $O_X$-module of sections of $E$ and
$\eE^{\vee}$ its dual.

Any cone $C=\spec\bigoplus_{i\geq0}S^i$ is canonically a closed
subcone of an abelian cone $A(C)=\spec\Sym S^1$, called the {\em
associated abelian cone }or the {\em abelian hull }of $C$. The abelian
hull is a vector bundle if and only if $S^1$ is locally free. Any
morphism of cones $\phi:C\rightarrow D$ induces a morphism
$A(\phi):A(C)\rightarrow A(D)$, extending $\phi$. Thus $A$ defines a
functor from cones to abelian cones called {\em abelianization}. Note
that $\phi$ is a closed immersion if and only if $A(\phi)$ is.

\begin{lem}
A cone $C$ over $X$ is a vector bundle if and only if it is smooth over $X$.
\end{lem}
\begin{pf} Let $C=\spec\bigoplus_{i\geq0}S^i$, and assume that $C\to X$ has
constant relative dimension $r$. Then $S^1=0^*\Omega_{C/X}$ is a rank $r$
vector
bundle. $C$ is a closed subcone of $A(C)=(S^1)^\vee$, hence by dimension
reasons
$C=A(C)$.
\end{pf}

If $E$ and $F$ are abelian cones over $X$, then any morphism of cones
$\phi:E\rightarrow F$ is a morphism of $X$-group schemes. If $E$ and
$F$ are vector bundles, then $\phi$ is a morphism of vector bundles.

\begin{example} If $X\rightarrow Y$ is a closed immersion with ideal sheaf
$\iI$, then
$$\bigoplus_{n\geq0}\iI^n/\iI^{n+1}$$ is a sheaf of $\O_X$-algebras and
$$C_{X/Y}=\spec
\bigoplus_{n\geq0}\iI^n/\iI^{n+1}$$ is a cone over $X$, called the {\em normal
cone }of $X$ in $Y$. The associated abelian cone $N_{X/Y}=\spec\Sym\iI/\iI^2$
is
also called the {\em normal sheaf }of
$X$ in $Y$.

More generally, any local immersion of Deligne-Mumford stacks has a normal cone
whose abelian hull is its normal sheaf (see \cite{vistoli}, definition 1.20).
\end{example}

\subsection{Exact Sequences of Cones}

\begin{defn} A sequence of cone morphisms $$
\ses E{i}C{}D$$
is {\em exact} if $E$ is a vector bundle and locally over $X$ there is a
morphism of cones $C\to E$ splitting $i$ and inducing an isomorphism $C\to
E\times D$.
\end{defn}
\begin{rmk} Given a short exact sequence $$
\ses {\fF'}{}\fF{}\eE$$
of coherent sheaves on $X$, with $\eE$ locally free, then $$
\ses {C(\eE)}{}{C(\fF')}{}{C(\fF)}$$
is exact, and conversely (see \cite{fulton}, Example 4.1.7).
\end{rmk}
\begin{lem} \label{ssmc}
Let $C\to D$ be a smooth, surjective morphism of cones, and let
$E=C\times_{D,0}X$; then the sequence $$\ses E{}C{}D$$ is exact.
\end{lem}
\begin{pf}
Write $C=\spec\bigoplus S^i$, $D=\spec\bigoplus S^{\prime i}$. We
start by proving that $$
\ses E{}{A(C)}{}{A(D)}$$is exact.

By base change we may assume $S^{\prime i}=0$ for $i\ge 2$. The cone
$E=\spec\Sym\eE$ is a vector bundle because it is smooth. On the other
hand, $E=\spec\bigoplus(S^i/S^{\prime 1}S^{i-1})$. As $C\to D$ is
smooth and surjective, $S^1\to S^{\prime 1}$ is injective. So we get
an exact sequence $$
\ses{S^1}{}{S^{\prime 1}}{}\eE.$$

To complete the proof, remark that $C\to A(C)\times_{A(D)}D$ is a
closed immersion, and both these schemes are smooth of the same
relative dimension over $C$.
\end{pf}

\subsection{$E$-Cones}

If $E$ is a vector bundle and $d:E\rightarrow C$ a morphism of cones, we say
that {\em $C$ is an
$E$-cone}, if $C$ is invariant under the action of $E$ on $A(C)$. We denote the
induced action of
$E$ on $C$ by
\begin{eqnarray*} E\times C & \longrightarrow & C \\ (\nu,\gamma) & \longmapsto
& d\nu+\gamma\quad.
\end{eqnarray*}
A {\em morphism }$\phi$ from an $E$-cone $C$ to an $F$-cone $D$ (or a
{\em morphism of vector bundle cones}) is a
commutative diagram of cones
\[\comdia{E}{d}{C}{\phi}{}{\phi}{F}{d}{D.}\] If
$\phi:(E,d,C)\rightarrow(F,d,D)$
and $\psi:(E,d,C)\rightarrow(F,d,D)$ are morphisms, we call them {\em
homotopic}, if there exists a morphism of cones $k:C\rightarrow F$, such that
\begin{enumerate}
\item $kd=\psi-\phi$,
\item $dk=\psi-\phi$.
\end{enumerate} Here the second condition is to be interpreted as saying that
$\phi+dk=\psi$. (More precisely, we say that $k$ is a {\em homotopy }from
$\phi$
to $\psi$.)

\begin{rmk} A sequence of cone morphisms with $E$ a vector bundle$$
\ses E{i}C{}D$$ is exact if and only if $C$ is an $E$-cone, $C\to D$ is
surjective, and the diagram $$\comdia{E\times
C}{\sigma}{C}{p}{}{\phi}{C}{\phi}{D}$$ is cartesian, where
$p$ is the projection and
$\sigma$ the action.
\end{rmk}

\begin{prop} \label{lcc}
Let $(C,0,\gamma)$ and $(D,0,\gamma)$ be algebraic $X$-spaces with
sections and $\aaa^1$-actions and let $\phi:C\rightarrow D$ be an
$\aaa^1$-equivariant $X$-morphism, which is smooth and surjective. Let
$E=C\times_{D,0}X$.  Then $C$ is an $E$-cone over $X$ if and only if
$D$ is a cone over $X$. Moreover, $C$ is abelian (a vector bundle) if
and only if $D$ is.
\end{prop}
\begin{pf}
Let us first assume that $C$ is an abelian cone, $C=\spec\Sym\fF$.
The morphism $E\rightarrow C$ gives rise to
$\fF\rightarrow\dual{\eE}$, where $\eE$ is the coherent $\O_X$-modules
of sections of $E$. Note that $\fF\rightarrow\dual{\eE}$ is an
epimorphism, since $E\rightarrow C$ is injective. Let $\gG$ be the
kernel, so that
\[\ses{\gG}{}{\fF}{}{\dual{\eE}}\] is a short exact sequence. Then
\[\ses{E}{}{C}{}{C(\gG)}\] is a short exact sequence of abelian cones over $X$,
so $D\cong C(\gG)$ and so $D$ is an abelian cone.

In general, $C\subset A(C)$ is defined by a homogeneous sheaf of ideals
$\iI\subset\Sym \sS^1$, where $\sS=\bigoplus\sS^i$ and $C=\spec\sS$. Let
$\fF=\sS^1$ and let $\gG$ as above be the kernel of
$\fF\rightarrow\dual{\eE}$. Let $\jJ=\iI\cap\Sym\gG$, which is a homogeneous
sheaf of ideals in
$\Sym\gG$, so $C'=\spec\Sym\gG/\jJ$ is a cone over $X$. By construction, $C'$
is
the scheme theoretic image of $C$ in $C(\gG)$. Hence $C'$ is the quotient of
$C$
by $E$ and so $C'\cong D$ and
$D$ is a cone.

Now for the converse. The claim is local in $X$. So since $D$ is affine over
$X$
we may assume that
$C=D\times E$ as $X$-schemes with $\aaa^1$-action. Then we are done.
\end{pf}

\subsection{Cone Stacks}

Let $X$ be, as above, a Deligne-Mumford stack over $k$. We need to
define the 2-category of algebraic stacks with $\aaa^1$-action over
$X$.

\begin{defn}
Let $\CC$ be an algebraic stack over $X$, together with a section
$0:X\to\CC$. An {\em $\aaa^1$-action }on $(\CC,0)$ is given by a
morphism of $X$-stacks $$\gamma:\aaa^1\times\CC\longrightarrow\CC$$ and
three 2-isomorphisms $\theta_1$, $\theta_0$ and $\theta_{\gamma}$
between the 1-morphisms in the following diagrams.
\begin{enumerate}
\item \[\comtri{\CC}{(1,\id)}{\aaa^1\times \CC}{\id}{\gamma}{\CC}\]
and $\theta_1:\id\to\gamma\comp(1,\id)$.
\item \[\comtri{\CC}{(0,\id)}{\aaa^1\times \CC}{0}{\gamma}{\CC}\]
and $\theta_0:0\to\gamma\comp(0,\id)$.
\item \[\comdia{\aaa^1\times\aaa^1\times \CC} {\id\times\gamma}
{\aaa^1\times \CC} {m\times\id}{}{\gamma} {\aaa^1\times \CC}{\gamma}{\CC}\]
and
$\theta_{\gamma}:\gamma\comp(m\times\id)\to\gamma\comp(\id\times\gamma)$.
\end{enumerate}
The 2-isomorphisms $\theta_1$, $\theta_0$ and $\theta_{\gamma}$ are
required to satisfy certain compatibilities which we leave to the
reader to make explicit (see also Section~1.4 in Expos\'e~XVIII of
\cite{sga4}, where a similar problem, the definition of Picard stacks,
is dealt with).

Let $(\CC,0,\gamma)$ and $(\DD,0,\gamma)$ be $X$-stacks with sections
and $\aaa^1$-actions. Then an {\em $\aaa^1$-equivariant morphism }
$\phi:\CC\to\DD$ is a triple $(\phi,\eta_0,\eta_{\gamma})$, where
$\phi:\CC\to\DD$ is a morphism of algebraic $X$-stacks and $\eta_0$
and $\eta_{\gamma}$ are 2-isomorphisms between the morphisms in the
following diagrams.
\begin{enumerate}
\item\begin{equation}\label{na1}
\comtri{X}{0}{\CC}{0}{\phi}{\DD}
\end{equation}
and $\eta_0:0\to\phi\comp 0$.
\item \begin{equation}\label{na2}
\comdia{\aaa^1\times\CC}{\id\times\phi}{\aaa^1\times\DD}
{\gamma}{}{\gamma} {\CC}{\phi}{\DD}
\end{equation}
and $\eta_{\gamma}:\phi\comp\gamma\to\gamma\comp(\id\times\phi)$.
\end{enumerate}
Again, the 2-isomorphisms have to satisfy certain compatibilities we
leave to the reader to spell out.

Finally, let $(\phi,\eta_0,\eta_{\gamma}):\CC\to\DD$ and
$(\psi,\eta'_0,\eta'_{\gamma}):\CC\to\DD$ be two $\aaa^1$-equivariant
morphisms. An {\em $\aaa^1$-equivariant isomorphism
}$\zeta:\phi\to\psi$ is a 2-isomorphism $\zeta:\phi\to\psi$ such that
the diagrams (notation compatible with (\ref{na1}) and (\ref{na2}))
\begin{enumerate}
\item \[\comtri{0}{\eta_0}{\phi\comp0}
{\eta'_0}{\zeta\comp0}{\psi\comp0} \]
\item \[\comdia{\phi\comp\gamma}{}{\gamma\comp(\id\times\phi)}
{\zeta\comp\gamma}{}{\gamma\comp(\id\times\zeta)}
{\psi\comp\gamma}{}{\gamma\comp(\id\times\psi)} \]
\end{enumerate}
commute.
\end{defn}

If $C$ is an $E$-cone, then since $E$ acts on $C$, we may form the
stack quotient of $C$ by $E$ over $X$, denoted $[C/E]$. For an
$X$-scheme $T$, the groupoid of sections of $[C/E]$ over $T$ is the
category of pairs $(P,f)$, where $P$ is an $E$-torsor (a principal
homogeneous $E$-bundle) over $T$ and $f:P\rightarrow C$ is an
$E$-equivariant morphism.

The $X$-stack $[C/E]$ comes with a section $0:X\to[C,E]$ and an
$\aaa^1$-action $\gamma:\aaa^1\times[C/E]\to[C/E]$. The section $0$ is
given by the pair $(E_T,0)$ over every $X$-scheme $T$; here $E_T$ is
the trivial $E$-bundle on $T$ and $0:E_T\to C$ is the vertex morphism.
The $\aaa^1$-action of $\alpha\in\aaa^1(T)=\O_T(T)$ on the category
$[C/E](T)$ is given by $\alpha\cdot(P,f)=(\alpha P,\alpha f)$, where
$\alpha P=P\times_{E,\alpha} E$ and $\alpha f:P\times_{E,\alpha} E\to
C$ is given by $[p,\nu]\mapsto \alpha f(p)+d(\nu)$.

If $\phi:(E,C)\to(F,D)$ is a morphism of vector bundle cones we get an
induced $\aaa^1$-equivariant morphism $\tilde{\phi}:[C/E]\to[D/F]$. A
homotopy $k:\phi\to\psi$ gives rise to an $\aaa^1$-equivariant
2-isomorphism $\tilde{k}:\tilde{\phi}\to\tilde{\psi}$ of
$\aaa^1$-equivariant morphism of stacks with $\aaa^1$-action. (See
Section~\ref{sohh} where these constructions are made explicit in a
similar case.)

\begin{lem}\label{lo2i}
Let $\phi,\psi:(E,C)\to(F,D)$ be morphisms and
$\zeta:\tilde{\phi}\to\tilde{\psi}$ an $\aaa^1$-equivariant
2-isomorphism between the associated $\aaa^1$-equivariant morphisms
$[C/E]\to [D/F]$. Then $\zeta=\tilde{k}$, for a unique
homotopy $k:\phi\to\psi$.
\end{lem}
\begin{pf}
We indicate how to construct $k:C\to F$. Given a section $c\in C(T)$
of $C$ over the $X$-scheme $T$, we consider the induced object
$(E_T,c)$ of $[C/E](T)$. The associated $F_T$-torsors
$E_T\times_{E_T,\phi^0}F_T$ and $E_T\times_{E_T,\psi^0}F_T$ are
trivial, so that $\phi(T)(E_T,c)$ is a section of $F$ over $T$. This
section we define to be $k(c)$.
\end{pf}

\begin{prop} \label{cics}
Let $C$ be an $E$-cone and $D$ an $F$-cone. Let
$\phi:(E,C)\rightarrow(F,D)$ be a morphism. If the diagram
\[\comdia{E}{}{C}{}{}{}{F}{}{D}\] is cartesian and $F\times C\rightarrow
D;(\mu,\gamma)\mapsto d\mu+\phi(\gamma)$ is surjective, then
$[C/E]\rightarrow[D/F]$ is an isomorphism of algebraic $X$-stacks with
$\aaa^1$-action.
\end{prop}
\begin{pf}
Similar to the proof of Proposition~\ref{tpifs} below.
\end{pf}

\begin{defn} \label{docs}
We call an algebraic stack $(\CC,0,\gamma)$ over $X$ with section and
$\aaa^1$-action a {\em cone stack}, if, locally with respect to the
\'etale topology on $X$, there exists a cone $C$ over $X$ and an
$\aaa^1$-equivariant morphism $C\to\CC$ that is smooth and surjective.

The morphism $C\rightarrow\CC$, or by abuse of language $C$, is called
a {\em local presentation }of $\CC$. The section $0:X\to\CC$ is called
the {\em vertex }of $\CC$.

Let $\CC$ and $\DD$ be cone stacks over $X$. A {\em morphism of cone
stacks } $\phi:\CC\rightarrow\DD$ is an $\aaa^1$-equivariant morphism
of algebraic $X$-stacks.

A {\em 2-isomorphism of cone stacks } is just an $\aaa^1$-equivariant
2-isomorphism.
\end{defn}

If $C\rightarrow\CC$ is a presentation of $\CC$, and
$E=C\times_{\CC,0}X$, then $C$ is an $E$-cone and $\CC\cong[C/E]$
as stacks with $\aaa^1$-action (use Lemma~\ref{ssmc} and
Proposition~\ref{lcc}).

If $\phi:\CC\to\DD$ is a morphism of cone stacks, then, locally with
respect to the \'etale topology on $X$, $\phi$ is
$\aaa^1$-equivariantly isomorphic to $[C/E]\to [D/F]$, where $E\to F$
is a morphism of vector bundles over $X$ and $C\to D$ is a
morphism from the $E$-cone $C$ to the $F$-cone $D$.

A 2-isomorphism of cone stacks $\zeta:\phi\to\psi$, where
$\phi,\psi:\CC\to\DD$, is locally over $X$ given by a homotopy of
morphisms of vector bundle cones. More precisely, one can find local
presentations $\CC\cong[C/E]$ and $\DD\cong[D/F]$ such that both
$\phi$ and $\psi$ are induced by morphisms of vector bundle cones
$\ol{\phi},\ol{\psi}:(E,C)\to(F,D)$ and under these identifications
$\zeta$ comes from a homotopy from $\ol{\phi}$ to $\ol{\psi}$. This
follows from Lemma~\ref{lo2i}.

\begin{rmk}
Let $\CC$ be a cone stack over $X$. By Proposition~\ref{lcc} the
fibered product over $\CC$ of any two local presentations is again a
local presentation. Moreover, if $\CC$ is a representable cone stack
over $X$, then $\CC$ is a cone. Every fibered product of cone stacks
is a cone stack.
\end{rmk}

\begin{examples}
All cones are cone stacks and all morphisms of cones are morphisms of
cone stacks. For a vector bundle $E$ on $X$, the classifying stack
$BE$ is a cone stack. Every homomorphism of vector bundles
$\phi:E\rightarrow F$ gives rise to a morphism of cone stacks.
\end{examples}

\begin{defn}
A cone stack $\CC$ over $X$ is called {\em abelian}, if, locally in
$X$, one can find presentations $C\rightarrow\CC$, where $C$ is an
abelian cone. A cone stack is a {\em vector bundle stack}, if one can
find such local presentations such that $C$ is a vector bundle. If
$\CC$ is abelian (a vector bundle stack), then for every local
presentation $C\rightarrow\CC$ the cone $C$ will be abelian (a vector
bundle).
\end{defn}

\begin{prop}
Every cone stack is a closed subcone stack of an abelian cone
stack. There exists a universal such abelian cone stack. It is called
the {\em associated abelian cone stack }or the {\em abelian hull}.
\end{prop}
\begin{pf}
Just glue the stacks obtained from the abelian hulls of local
presentations.
\end{pf}

\begin{defn} \label{dics}
Let $\EE$ be a vector bundle stack and $\EE\to\CC$ a morphism of cone
stacks. We say that $\CC$ is an {\em $\EE$-cone stack}, if $\EE\to\CC$
is locally isomorphic (as a morphism of cone stacks, i.e.\
$\aaa^1$-equivariantly) to the morphism $[E_1/E_0]\to[C/F]$ coming from
a commutative diagram
\[\comdia{E_0}{}{F}{}{}{}{E_1}{}{C,}\]
where $C$ is both an $E_1$- and an $F$-cone.
\end{defn}

If $\CC$ is an $\EE$-cone stack, then there exists a natural morphism
$\EE\times\CC\to\CC$ coming from the action $E_1\times C\to C$ in a
local presentation of $\EE\to\CC$ as above.  We call
$\EE\times\CC\to\CC$ the {\em action } of $\EE$ on $\CC$.

\begin{defn} \label{dsescs}
Let $\EE\to\CC\to\DD$ be a sequence of morphisms of cone stacks, where
$\CC$ is an $\EE$-cone stack. If
\begin{enumerate}
\item $\CC\to\DD$ is a smooth epimorphism,
\item the diagram
\[\comdia{\EE\times\CC}{\sigma}{\CC}{p}{}{}{\CC}{}{\DD}\]
(where $p$ is the projection and $\sigma$ the action) is cartesian,
\end{enumerate}
we call $0\to\EE\to\CC\to\DD\to 0$ a {\em short exact sequence }of
cone stacks.  Note that this is equivalent to $\CC$ being locally
isomorphic to $\EE\times \DD$.\end{defn}
\begin{prop}\label{csescs}
The sequence $\EE\to\CC\to\DD$ of morphisms of cone stacks is exact if
and only if locally in $X$ there exist commutative diagrams
\[\begin{array}{ccccccccc}
0 & \longrightarrow & E_0 & \longrightarrow & F & \longrightarrow & G
& \longrightarrow & 0 \\
& &\downarrow & & \downarrow  & & \downarrow & & \\
0 &\longrightarrow & E_1 & \longrightarrow & C & \longrightarrow & D &
\longrightarrow& 0,
\end{array}\]
where the top row is a short exact sequence of vector bundles and the
bottom row is a short exact sequence of cones, such that
$\EE\to\CC\to\DD$ is isomorphic to $[E_1/E_0]\to[C/F]\to[D/G]$.
\end{prop}
\begin{pf}
The statement is local on $X$. To prove the only if part we can assume
$\CC=\EE\times \DD$, and then it's trivial. To prove the if part, note
that both short exact sequences are locally split.
\end{pf}

\section{Stacks of the Form $h^1/h^0$} \label{sohh}

\subsection{The General Theory}

We shall review here some aspects of the theory of Picard stacks
developed by Deligne in Section~1.4 of Expos\'e~XVIII in \cite{sga4}.
For the precise definition of Picard stack see [ibid.]. Roughly
speaking, a Picard stack is a stack together with an `addition'
operation, that is both associative and commutative. An example would
be the stack of torsors under a commutative group sheaf.

Let $X$ be a topos and $d:E^0\rightarrow E^1$ a homomorphism of
abelian sheaves on $X$, which we shall consider as a complex of
abelian sheaves on $X$. Via $d$, the abelian sheaf $E^0$ acts on $E^1$
and we may consider the stack-theoretic quotient of this action,
denoted
\[h^1/h^0(E\com)=[E^1/E^0],\]
which is a Picard stack on $X$. (See also [ibid.] 1.4.11, where
$h^1/h^0(E\com)$ is denoted $\mbox{ch}(E\com)$.)  For an object
$U\in\ob X$ the groupoid $h^1/h^0(E\com)(U)$ of sections of
$h^1/h^0(E\com)$ over $U$ is the category of pairs $(P,f)$, where $P$
is an $E^0$-torsor (principal homogeneous $E^0$-bundle) over $U$ and
$f:P\rightarrow E^1\resto U$ is an $E^0$-equivariant morphism of
sheaves on $U$.

Now if $d:F^0\rightarrow F^1$ is another homomorphism of abelian
sheaves on $X$ and $\phi:E\com\rightarrow F\com$ is a homomorphism of
homomorphisms (or in other words a homomorphism of complexes), then we
get an induced morphism of Picard stacks (an additive morphism in the
terminology of [ibid.])
\[h^1/h^0(\phi):h^1/h^0(E\com)\longrightarrow h^1/h^0(F\com).\] For an object
$U\in\ob X$ the functor $h^1/h^0(\phi)(U)$ maps the pair $(P,f)$ to the pair
$(P\times_{E^0,\phi^0}F^0,\phi^1(f))$, where $\phi^1(f)$ denotes the map
\begin{eqnarray*}
\phi^1(f):P\times_{E^0}F^0 & \longrightarrow & F^1 \\ {[p,\nu]} & \longmapsto &
\phi^1(f(p))+d(\nu).
\end{eqnarray*}

Now, if $\psi:E\com\rightarrow F\com$ is another homomorphism of complexes and
$k:\phi\rightarrow\psi$ is a homotopy, i.e.\ a homomorphism of abelian sheaves
$k:E^1\rightarrow F^0$, such that
\begin{enumerate}
\item $kd=\psi^0-\phi^0$,
\item $dk=\psi^1-\phi^1$,
\end{enumerate} then we get an induced isomorphism
$\theta:h^1/h^0(\phi)\rightarrow h^1/h^0(\psi)$ of morphisms of Picard
stacks from $h^1/h^0(E\com)$ to $h^1/h^0(F\com)$. If $U\in\ob X$ is an
object, then $\theta(U)$ is a natural transformation of functors from
$h^1/h^0(\phi)(U)$ to $h^1/h^0(\psi)(U)$. For an object $(P,f)$ of
$h^1/h^0(E\com)(U)$ the morphism $\theta(U)(P,f)$ is a morphism from
$h^1/h^0(\phi)(U)(P,f)$ to $h^1/h^0(\psi)(U)(P,f)$ in the category
$h^1/h^0(F\com)(U)$. In fact, $\theta(U)(P,f)$ is the isomorphism of
$F^0\resto U$-torsors
\begin{eqnarray}
\theta(U)(P,f):P\times_{E^0,\phi^0}F^0 & \longrightarrow &
P\times_{E^0,\psi^0}F^0\\ {[p,\nu]} & \longmapsto & [p,kf(p)+\nu],\nonumber
\end{eqnarray} such that the diagram of $F^0\resto U$-sheaves
\[\inversecomtri{P\times_{E^0,\phi^0}F^0}{\theta(U)(P,f)}{\phi^1(f)
}{P\times_{E^0,\psi^0}F^0}
{\psi^1(f)}{F^1}\] commutes.

\begin{prop} \label{tpifs}
Let $\phi:E\com\rightarrow F\com$ be a homomorphism of homomorphisms
of abelian sheaves on $X$, as above. If $\phi$ induces isomorphisms on
kernels and cokernels (i.e.\ if $\phi$ is a quasi-isomorphism), then
$h^1/h^0(\phi):h^1/h^0(E\com)\rightarrow h^1/h^0(F\com)$ is an
isomorphism of Picard stacks over $X$.
\end{prop}
\begin{pf}
First let us treat the case that $\phi$ is a homotopy equivalence.
Then, in fact, any homotopy inverse of $\phi$ will provide an inverse
to $h^1/h^0(\phi)$, by the above remarks.

As a second case, let us assume that $\phi\com:E\com\rightarrow F\com$ is an
epimorphism (i.e.\
$\phi^0$ and $\phi^1$ are epimorphisms). In this case $E^1\rightarrow[F^1/F^0]$
is an epimorphism, so for $[E^1/E^0]$ to be isomorphic to $[F^1/F^0]$, it is
necessary and sufficient that
\[\comdia{E^0\times E^1}{d+\id}{E^1}{\pr}{}{}{E^1}{}{[F^1/F^0]}\] be cartesian.
This quickly reduces to proving that
\[\comdia{E^1\times E^0}{}{E^1}{}{}{}{E^1\times F^0}{}{F^1}\] is cartesian,
which, in turn, is equivalent to
\[\comdia{E^0}{}{E^1}{}{}{}{F^0}{}{F^1}\] being cartesian, which is a
consequence of the assumptions.

Finally, let us note that a general $\phi$ factors as a homotopy equivalence
followed by an epimorphism. To see this consider $E\com\oplus F^0$, which is
homotopy equivalent to $E\com$. Define a homomorphism $\psi:E\com\oplus
F^0\rightarrow F\com$ by $\psi^0(\nu,\mu)=\phi^0(\nu)+\mu$ and
$\psi^1(\chi,\mu)=\phi^1(\chi)+d(\mu)$. Then $\psi$ is surjective and
$\phi=\psi\comp i$, where
$i:E\com\rightarrow E\com\oplus F^0$ is given by $i=\id\oplus 0$.
\end{pf}

If $E\com$ is a complex of arbitrary length of abelian sheaves on $X$, let
\begin{eqnarray*} Z^i(E\com) & = & \ker(E^i\rightarrow E^{i+1}) \\ C^i(E\com) &
= & \cok(E^{i-1}\rightarrow E^i).
\end{eqnarray*} The complex $E\com$ induces a homomorphism
\[\tau_{[0,1]} E\com=[C^0(E\com)\rightarrow Z^1(E\com)]\] and we let
$h^1/h^0(E\com)=h^1/h^0(\tau_{[0,1]}E\com)$.

Now let $\O_X$ be a sheaf of rings on $X$ and $C(\O_X)$, $K(\O_X)$ and
$D(\O_X)$
the category of complexes of $\O_X$-modules, the category of complexes of
$\O_X$-modules up to homotopy and the derived category of the category
$\Mod(\O_X)$ of $\O_X$-modules, respectively. Let
$\phi:E\com\rightarrow F\com$ be a morphism in $D(\O_X)$. Let
\[\begin{array}{ccc} H\com & \stackrel{\psi}{\longrightarrow} & F\com \\
\ldiag{\alpha} & & \\ E\com & &
\end{array}\] be a diagram in $C(\O_X)$ giving rise to $\phi$, where $\alpha$
is a quasi-isomorphism. We get an induced diagram of Picard stacks
\[\begin{array}{ccc} h^1/h^0(H\com) & \stackrel{h^1/h^0(\psi)}{\longrightarrow}
& h^1/h^0(F\com) \\
\ldiag{h^1/h^0(\alpha)} & & \\ h^1/h^0(E\com), & &
\end{array}\] where $h^1/h^0(\alpha)$ is an isomorphism by
Proposition~\ref{tpifs}. Choosing an inverse of
$h^1/h^0(\alpha)$ induces a morphism
\[h^1/h^0(E\com)\longrightarrow h^1/h^0(F\com).\] One checks that different
choices of $(\alpha, H\com,\psi)$ and $h^1/h^0(\alpha)^{-1}$ give rise
to isomorphic morphisms $h^1/h^0(E\com)\rightarrow
h^1/h^0(F\com)$. This proves in particular that if $E\com$ and $F\com$
are isomorphic in $D(\O_X)$, then the Picard $X$-stacks
$h^1/h^0(E\com)$ and $h^1/h^0(F\com)$ are isomorphic.

\begin{example}
If $d:E^0\rightarrow E^1$ is a monomorphism then
$h^1/h^0(E\com)=\cok(d)$ is a sheaf over $X$.

If $d:E^0\rightarrow E^1$ is an epimorphism then
$h^1/h^0(E\com)=B\ker(d)$ is a gerbe over $X$.
\end{example}

\begin{lem} \label{lmh}
1. Let $\phi,\psi:E\com\to F\com$ be two morphisms in $D(\O_X)$. Then,
if for some choice of $h^1/h^0(\phi)$ and $h^1/h^0(\psi)$ we have
$h^1/h^0(\phi)\cong h^1/h^0(\psi)$ as morphisms of Picard stacks, then
$\phi=\psi$.

2. Let $0(E,F)$ be the zero morphism $0(E,F):h^1/h^0(E\com)\to
h^1/h^0(F\com)$. Then $\Aut(0(E,F))=\Hom^{-1}_{D(\O_X)}(E\com,F\com)$.
\end{lem}
\begin{pf}
These are similar to Lemma~\ref{lo2i}. See also [ibid.].
\end{pf}

\subsection{Application to Schemes}

Let $X$ be a Deligne-Mumford stack. Consider the morphism of topoi
\[v:X_\fl\longrightarrow X_\et.\] The functor $v\lst$ restricts a sheaf on the
big fppf-site to the small \'etale site and its left adjoint $v^{-1}$ extends
the embedding of the \'etale site into the flat site.

Let $\O_{X_\fl}$ and $\O_{X_\et}$ denote the sheaves of rings induced by $\O_X$
on $X_\fl$ and
$X_\et$, respectively. There is a canonical morphism of sheaves of rings
$v^{-1}\O_{X_\et}\rightarrow\O_{X_\fl}$, so that we have a morphism of ringed
topoi
$$v:(X_\fl,\O_{X_\fl})\rightarrow (X_{\et},\O_{X_\et}).$$ The induced functor
{}from $\Mod(\O_{X_\et})$ to $\Mod(\O_{X_\fl})$ will be denoted by $v\upst$:
\[v\upst(M)=v^{-1}M\otimes_{v^{-1}\O_{X_\et}}\O_{X_\fl}.\] Since
$\Mod(\O_{X_\et})$ has enough flat modules we may derive the right exact
functor
$v\upst$ to get the functor $Lv\upst:D^-(\O_{X_\et})\rightarrow
D^-(\O_{X_\fl})$. To abbreviate notation, we write $M_\fl\com=Lv\upst M\com$
for
$M\com\in\ob D^-(\O_{X_\et})$.

We shall also need to consider the functor
\[R\sheafhom(\argument,\O_{X_\fl}):D^-(\O_{X_\fl})\longrightarrow
D^+(\O_{X_\fl}).\] It is defined using an injective resolution
$\O_{X_\fl}\iso\iI\com$ of $\O_{X_\fl}$, i.e.\
\[R\sheafhom(M\com,\O_{X_\fl})=\tot\sheafhom(M\com,\iI\com),\] but if $M\com$
happens to have a projective resolution $\pP\com\iso M\com$, then we have
\[R\sheafhom(M\com,\O_{X_\fl})\cong\sheafhom(\pP\com,\O_{X\fl}).\] We shall
abbreviate notation by writing
\[\dual{M\com}=R\sheafhom(M\com,\O_{X_\fl}).\]

We will be interested in the stack $h^1/h^0(\dual{(M\com_\fl)})$ associated to
an object
$M\com\in\ob D^-(\O_{X_\et})$. Note that for such $M\com\in\ob D^-(\O_{X_\et})$
we have
\[h^1/h^0(\dual{(M\com_\fl)})\cong h^1/h^0(\dual{(\tau_{\geq-1}M\com_\fl)}).\]

\begin{defn} \label{dost} We say that an object $L\com$ of
$D(\O_{X_\et})$ {\em satisfies Condition~}($\star$) if
\begin{enumerate}
\item $h^i(L\com)=0$ for all $i>0$,
\item $h^i(L\com)$ is coherent, for $i=0,-1$.
\end{enumerate}
\end{defn}

\begin{prop} \label{has} Let $L\com\in\ob D(\O_{X_\et})$ satisfy
Condition~{\rm (}$\star${\rm )}. Then the $X$-stack
$h^1/h^0(\dual{(L\com_\fl)})$ is an algebraic $X$-stack, in fact an abelian
cone
stack over $X$. Moreover, if $L\com$ is of perfect amplitude contained in
$[-1,0]$, then
$h^1/h^0(\dual{(L\com_\fl)})$ is a vector bundle stack.
\end{prop}
\begin{pf} The claim is local in $X$ (with respect to the \'etale topology), so
we may assume that $L\com$ has a free resolution, or that $L\com$ itself
consists of free $\O_X$-modules. We may also assume that
$L^i=0$, for all $i>0$ and that $L^0$ and $L^{-1}$ have finite rank. Then
$L\com_\fl$ is given by
$L\com$ itself, since a free sheaf is flat, and $\dual{(L\com_\fl)}$ is given
by
$\dual{L}\com$, taking duals component-wise, since a free module is projective.
Thus
\[h^1/h^0(\dual{(L\com_\fl)})=[Z^1(\dual{L}\com)/\dual{L}^0],\] which is the
cone stack given by the homomorphism of abelian cones $\dual{L}^0\rightarrow
Z^1(\dual{L}\com)=C(C^{-1}(L\com))$.

If $L\com$ is of perfect amplitude contained in $[-1,0]$, then we may assume
that in addition to the above assumptions $L^i=0$, for all $i\leq-2$. Then
$Z^1(\dual{L}\com)=\dual{L}^1$ is a vector bundle.
\end{pf}

So if $\phi:E\com\rightarrow L\com$ is a homomorphism in $D(\O_{X_\et})$, where
$E\com$ and
$L\com$ satisfy ($\star$), then we get an induced morphism of algebraic stacks
\[\dual{\phi}:h^1/h^0(\dual{(L\com_\fl)})\longrightarrow
h^1/h^0(\dual{(E\com_\fl)}).\]

\begin{prop}\label{rmcs}
The morphism $\dual{\phi}$ is a morphism of abelian cone stacks.
Moreover, $h^0(\phi)$ is surjective, if and only if $\dual{\phi}$ is
representable.
\end{prop}
\begin{pf}
The fact that $\dual{\phi}$ is a morphism of abelian cone stacks is
immediate from the definition. The second question is local in $X$, so
we may assume that $E\com$ and $L\com$ are complexes of free
$\O_X$-modules and that $E^i=L^i=0$, for $i>0$, and that $L^0$,
$L^{-1}$, $E^0$ and $E^{-1}$ are of finite rank. Consider the
commutative diagram
\[\comdia{C^{-1}(E\com)}{}{E^0}{}{}{}{C^{-1}(L\com)}{}{L^0}\] of coherent
sheaves on $X$. Let
$F$ be the fibered product
\[\comdia{F}{}{E^0}{}{}{}{C^{-1}(L\com)}{}{L^0.}\] The fact that $h^0(\phi)$ is
surjective, is equivalent to saying that the sequence
\[\ses{F}{}{E^0\oplus C^{-1}(L\com)}{}{L^0}\]
is exact. Since $L^0$ is free, we get an
induced exact sequence of cones
\[\ses{\dual{L}^0}{}{\dual{E}^0\oplus Z^1(\dual{L}\com)}{}{C(F)}.\] Hence
by Proposition~\ref{cics} we have
\[[Z^1(\dual{L}\com)/\dual{L}^0]\cong[C(F)/\dual{E}^0].\]
In particular the diagram
\[\comdia{C(F)}{}{Z^1(\dual{E}\com)}{}{}{}{h^1/h^0(\dual{(L\com_\fl)})}{}
{h^1/h^0(\dual{(E\com_\fl)})}\] is cartesian, hence $\dual{\phi}$ is
representable.

For the converse, note that $\dual{\phi}$ representable implies that
$\dual{L}^0\to\dual{E}^0\times Z^1(\dual{L}\com)$ is a closed
immersion, which implies that $E^0\oplus C^{-1}(L\com)\to L^0$ is an
epimorphism.
\end{pf}

\begin{prop}\label{isoiso}
The morphism $\dual{\phi}$ is a closed immersion if and only if
$h^0(\phi)$ is an isomorphism and $h^{-1}(\phi)$ is surjective.
Moreover, $\dual{\phi}$ is an isomorphism if and only if $h^0(\phi)$
and $h^{-1}(\phi)$ are.
\end{prop} \begin{pf}
Following the previous argument, $\dual{\phi}$ is a closed immersion
if and only if $C(F)\to Z^1(\dual{E}\com)$ is. This is equivalent to
$C^{-1}(E\com)\to F$ being surjective. A simple diagram chase shows
that this is equivalent to $h^0(\phi)$ being an isomorphism and
$h^{-1}(\phi)$ being surjective. The `moreover' follows similarly.
\end{pf}

\begin{prop} \label{dtscs}
Let
\[\dt{E\com}{}{F\com}{}{G\com}{}\]
be a distinguished triangle in $D(\O_{X_\et})$, where $E\com$ and
$F\com$ satisfy ($\star$) and $G\com$ is of perfect amplitude
contained in $[-1,0]$. Then the induced sequence
\[h^1/h^0(\dual{G})\longrightarrow h^1/h^0(\dual{F})\longrightarrow
h^1/h^0(\dual{E})\]
is a short exact sequence of abelian cone stacks over $X$.
\end{prop}
\begin{pf} The question is local, so assume that
$E^i$ and $F^i$ are $0$ for $i>0$ and vector bundles for $i=0,-1$,
and that
$G^i=F^i\oplus E^{i+1}$. We have to prove that $$
\ses {[Z^1(\dual{G})/\dual{G}^0]} {} {[Z^1(\dual{F})/\dual{F}^0]}{}
{[Z^1(\dual{E})/\dual{E}^0]}$$
is a short exact sequence of cone stacks.
By Proposition \ref{csescs}, it is enough to prove that the exact sequence of
sheaves
$$
\ses {C^{-1}(E\com)} {} {C^{-1}(F\com)\oplus E^0}{} {C^{-1}(G\com)}$$
is exact. This is then a straightforward verification.
\end{pf}

\section{The Intrinsic Normal Cone} \label{stinc}

\subsection{Normal Cones}

Normal cones have the following functorial property. Consider a commutative
diagram of (arbitrary) algebraic $k$-stacks
\begin{equation} \label{cdci}
\comdia{X'}{j}{Y'}{u}{}{v}{X}{i}{Y,}
\end{equation} where $i$ and $j$ are local immersions. Then there is
a natural morphism of cones over $X'$
\[\alpha:C_{X'/Y'}\longrightarrow u\upst C_{X/Y}.\] If (\ref{cdci}) is
cartesian, then $\alpha$ is a closed immersion. If, moreover, $v$ is
flat, then $\alpha$ is an isomorphism.

\begin{prop} \label{lonc} Consider a commutative diagram of Deligne-Mumford
stacks
\[\comtri{X}{i'}{Y'}{i}{f}{Y,}\] where $i$ and $i'$ are local immersions
and
$f$ is smooth. Then the sequence of morphisms of cones over $X$
\begin{equation} \label{sesc} {i'}\upst
T_{Y'/Y}\stackrel{\beta}{\longrightarrow}C_{X/Y'}\stackrel{
\alpha}{\longrightarrow}C_{X/Y},
\end{equation}
where the maps $\alpha$ and $\beta$ are the natural ones, is exact.
\end{prop}
\begin{pf} The question is local, so we can assume that $X$, $Y$ and $Y'$ are
schemes and that $i'$ and $i$ are immersions. This is then Example 4.2.6 in
\cite{fulton}.
\end{pf}

\begin{lem}\label{nci} %normal cone invariant
Let
\[U \stackrel{f}{\longrightarrow}  M \] be a local immersion of affine
$k$-schemes of finite type, where $M$ is smooth over $k$. Then the normal cone
$C_{U/M}\hookrightarrow N_{U/M}$ is invariant under the action of $f\upst T_M$
on
$N_{U/M}$. In other words, $C_{U/M}$ is an $f\upst T_M$-cone.
\end{lem}
\begin{pf}
Let $p_i:M\times M\rightarrow M$, $i=1,2$, be the two projections.
Each one gives rise to a commutative diagram
\[\comtri{U}{\Delta f}{M\times M}{f}{p_i}{M,}\]
and hence to an exact sequence
\[\ses{f\upst T_M}{j_i}{N_{U/M\times M}}{{p_i}\lst}{N_{U/M}}\]
of abelian cones on $U$.

The diagonal gives rise to the commutative diagram
\[\comtri{U}{f}{M}{\Delta f}{\Delta}{M\times M}\]
and hence to a homomorphism
\[N_{U/M}\stackrel{s}{\longrightarrow} N_{U/M\times M}\]
of abelian cones on $U$.

Now $s$ is a section of both ${p_1}\lst$ and ${p_2}\lst$. Using
$(j_1,{p_1}\lst)$ we make the identification
\begin{equation} \label{nci1}
N_{U/M \times M}=f\upst T_M\times N_{U/M}.
\end{equation}
Then ${p_2}\lst$ is identified with the action of $f\upst T_M$ on
$N_{U/M}$. Since the same functorialities of normal sheaves used so
far are enjoyed by normal cones, we get that under the
identification~(\ref{nci1}) the subcone $C_{U/M\times M}\subset
N_{U/M\times M}$ corresponds to $f\upst T_M\times C_{U/M}$ and the
action ${p_2}\lst:f\upst T_M\times N_{U/M}\rightarrow N_{U/M}$
restricts to ${p_2}\lst:f\upst T_M\times C_{U/M}\rightarrow C_{U/M}$.
\end{pf}

The following is not used until Section~\ref{otafc}.

Consider the diagram~(\ref{cdci}), assume it is cartesian and assume
that $v$ is a regular local immersion. Assume also that $Y$ is smooth
of constant dimension.  Let $C=C_{X/Y}$ and $N=N_{Y'/Y}$. Then we get
an induced cartesian diagram
\begin{equation} \label{bdnc}
\begin{array}{ccccc}
N\times_YC & \longrightarrow & u\upst C & \longrightarrow & C \\
\ldiag{} & & \ldiag{} & & \rdiag{} \\
j\upst N & \longrightarrow & X' & \stackrel{u}{\longrightarrow} & X \\
\ldiag{} & & \ldiag{j} & & \rdiag{i} \\
N & \stackrel{\rho}{\longrightarrow} & Y' &
\stackrel{v}{\longrightarrow} & Y.
\end{array}
\end{equation}
If $Y$ is a scheme, Vistoli constructed in \cite{vistoli} a canonical
rational equivalence $\beta(Y',X)\in W\lst(N\times_YC)$ such that
\[\del\beta(Y',X)=[C_{u\upst C/C}]-[\rho\upst C_{X'/Y'}].\]

\begin{note}
Let $0:u\upst C\rightarrow N\times_YC$ be the zero section. Then
$$0\upsh[C_{u\upst C/C}]=v\upsh[C]\in A\lst(u\upst C),$$ by the
definition of $v\upsh$. On the other hand, $$0\upsh[\rho\upst
C_{X'/Y'}]=0\upsh\rho\upsh[C_{X'/Y'}]=[C_{X'/Y'}]\in A\lst(u\upst
C).$$ So the existence of Vistoli's rational equivalence implies that
$$v\upsh[C]=[C_{X'/Y'}].$$
\end{note}

\begin{prop}  \label{vrqsb}
Vistoli's rational equivalence commutes with any smooth base change
$\phi:Y_1\rightarrow Y$. More precisely, if we denote by a subscript
$(\argument)_1$ the base change via $\phi$ of any object in
(\ref{bdnc}), then
$$\phi\upst\beta(Y',X)=\beta(Y'_1,X_1)\in W\lst(N_1\times_{Y_1}C_1).$$
\end{prop}
\begin{pf}
If $\phi$ is \'etale, this is Lemma 4.6(ii) in \cite{vistoli}. Vistoli's proof
is based on the fact that the following commute with \'etale base change:
blowing up a scheme along a closed subscheme; normalization; order of
a Cartier divisor along an irreducible Weil divisor on a reduced,
equidimensional
scheme. But all these operations do in fact commute with smooth base change.
\end{pf}

A first consequence of this proposition is that we may drop the
assumption that $Y$ be a scheme. We get $\beta(Y',X)\in
W\lst(N\times_YC)$ for any situation~(\ref{bdnc}). The consequence
$v\upsh[C]=[C_{X'/Y'}]$ holds if $Y$ (and hence all other stacks in
(\ref{bdnc})) is of Deligne-Mumford type.

Now let us assume that $i:X\rightarrow Y$ factors as
\[\comtri{X}{\tilde{\imath}}{\tilde{Y}}{i}{\pi}{Y,}\]
where $\tilde{\imath}$ is another local immersion and $\pi$ is of
relative Deligne-Mumford type (i.e.\ has unramified diagonal) and is
smooth of constant fiber dimension.  Then we construct the cartesian
diagram
\[\comdia{\tilde{Y}'}{\tilde{v}}{\tilde{Y}}{}{}{\pi}{Y'}{v}{Y}\]
and over
\[\comdia{X'}{u}{X}{\tilde{\jmath}}{}{\tilde{\imath}}{\tilde{Y}'}{\tilde{u}}
{\tilde{Y}}\]
we construct the analogue of (\ref{bdnc}):
\begin{equation} \label{bdnct}
\begin{array}{ccccc}
N\times_Y\tilde{C} & \longrightarrow & u\upst \tilde{C} &
\longrightarrow & \tilde{C} \\
\ldiag{} & & \ldiag{} & & \rdiag{} \\
j\upst N & \longrightarrow & X' & \stackrel{u}{\longrightarrow} & X \\
\ldiag{} & & \ldiag{\tilde{\jmath}} & & \rdiag{\tilde{\imath}} \\
\pi\upst N & \stackrel{\tilde{\rho}}{\longrightarrow} & \tilde{Y}' &
\stackrel{\tilde{v}}{\longrightarrow} & \tilde{Y},
\end{array}
\end{equation}
i.e.\ $\tilde{C}=C_{X/\tilde{Y}}$. Diagrams~(\ref{bdnc})
and~(\ref{bdnct}) may be fused into one large diagram
\begin{equation} \label{bdnctt}
\begin{array}{ccccc}
N\times_Y\tilde{C} & \longrightarrow & u\upst \tilde{C} &
\longrightarrow & \tilde{C} \\
\ldiag{} & & \ldiag{} & & \rdiag{\alpha} \\
N\times_YC & \longrightarrow & u\upst C & \longrightarrow & C \\
\ldiag{} & & \ldiag{} & & \rdiag{} \\
j\upst N & \longrightarrow & X' & \stackrel{u}{\longrightarrow} & X \\
\ldiag{} & & \ldiag{\tilde{\jmath}} & & \rdiag{\tilde{\imath}} \\
\pi\upst N & \stackrel{\tilde{\rho}}{\longrightarrow} & \tilde{Y}' &
\stackrel{\tilde{v}}{\longrightarrow} & \tilde{Y} \\
\ldiag{} & & \ldiag{} & & \rdiag{\pi} \\
N & \stackrel{\rho}{\longrightarrow} & Y' &
\stackrel{v}{\longrightarrow} & Y.
\end{array}
\end{equation}
By Proposition~\ref{lonc} the morphism $\tilde{C}\rightarrow C$ is a
$T_{\tilde{Y}/Y}\times_{\tilde{Y}}C$-bundle.

\begin{prop} \label{vcgt}
We have $\alpha\upst(\beta(Y',X))=\beta(\tilde{Y}',X)$ in
$W\lst(N\times_Y\tilde{C})$.
\end{prop}
\begin{pf}
By the compatibilities of $\beta$ proved in \cite{vistoli} we reduce
to the case that $\tilde{Y}=\aaa^n_Y$, $\pi:\aaa^n_Y\rightarrow Y$ is a
relative affine $n$-space and $\tilde{\imath}:Y\to\aaa^n_Y$ is the zero
section. Then one checks that Vistoli's construction commutes with
$\pi$.
\end{pf}

\begin{prop} \label{vrqi}
In the situation of Diagram~(\ref{bdnc}) assume that $Y$ is of
Deligne-Mumford type. Vistoli's rational equivalence $\beta(Y',X)\in
W\lst(N\times_YC)$ is invariant under the natural action of $j\upst
N\times_Y T_Y$ on $N\times_YC$.
\end{prop}
\begin{pf}
The vector bundle $i\upst T_Y$ acts on the $X$-cone $C$ by
Lemma~\ref{nci}. Pulling back from $X$ to $j\upst N$ gives the natural
action of $j\upst N\times_Y T_Y$ on $N\times_YC$. Using the
construction of the proof of Lemma~\ref{nci} the claim follows from
Proposition~\ref{vcgt} applied to $\tilde{Y}=Y\times Y$ and
$\tilde{\imath}=\Delta\comp i:X\to Y\times Y$.
\end{pf}

\subsection{The Intrinsic Normal Cone}

Let $X$ be a Deligne-Mumford stack, locally of finite type over $k$. Let
$L_X\com$ be the cotangent complex of $X$ relative to $k$. Then $L_X\com\in\ob
D(\O_{X\et})$ and $L_X\com$ satisfies ($\star$).

\begin{defn}
We denote the algebraic stack
$h^1/h^0(\dual{((L_X\com)_\fl)})$ by $\NN_X$ and call it the {\em
intrinsic normal sheaf }of $X$.
\end{defn}

We shall now construct the intrinsic normal cone as a closed subcone
stack of $\NN_X$.

\begin{defn}  A {\em local embedding }of $X$ is a diagram
\[\begin{array}{ccc} U & \stackrel{f}{\longrightarrow} & M \\
\ldiag{i} & & \\ X & &\quad,
\end{array}\] where
\begin{enumerate}
\item $U$ is an affine $k$-scheme of finite type,
\item $i:U\rightarrow X$ is an \'etale morphism,
\item $M$ is a smooth affine $k$-scheme of finite type,
\item $f:U\rightarrow M$ is a local immersion.
\end{enumerate} By abuse of language we call the pair $(U,M)$ a local embedding
of $X$.

A morphism of local embeddings $\phi:(U',M')\rightarrow(U,M)$ is a pair of
morphisms
$\phi_U:U'\rightarrow U$ and $\phi_M:M'\rightarrow M$ such that
\begin{enumerate}
\item $\phi_U$ is an \'etale $X$-morphism,
\item $\phi_M$ is a smooth morphism such that
\[\comdia{U'}{f'}{M'}{\phi_U}{}{\phi_M}{U}{f}{M}\] commutes.
\end{enumerate}
\end{defn}

If $(U',M')$ and $(U,M)$ are local embeddings of $X$, then
$(U'\times_X U,M'\times M)$ is naturally a local embedding of $X$
which we call the {\em product }of $(U',M')$ and $(U,M)$, even though
it may not be the direct product of $(U',M')$ and $(U,M)$ in the
category of local embeddings of $X$.

Let
\[\begin{array}{ccc} U & \stackrel{f}{\longrightarrow} & M \\
\ldiag{i} & & \\ X & &\quad
\end{array}\]
be a local embedding of $X$. Let $I/I^2$ be the conormal sheaf of $U$
in $M$. There is a natural homomorphism of coherent $\O_U$-modules
$I/I^2\rightarrow f\upst\Omega_M$. Moreover, there exists a natural
homomorphism
\[\phi:L_X\com\resto U\longrightarrow[I/I^2\rightarrow f\upst\Omega_M] \] in
$D(\O_{U_\et})$, where we think of $[I/I^2\rightarrow f\upst\Omega_M]$
as a complex concentrated in degrees $-1$ and $0$.  Moreover, $\phi$
induces an isomorphism on $h^{-1}$ and $h^0$ (see \cite{Ill},
Chapitre~III, Corollaire~3.1.3). Hence by Proposition~\ref{isoiso} we
get an induced isomorphism of cone stacks
\[\dual{\phi}:[N_{U/M}/f\upst T_M]\longrightarrow i\upst \NN_X,\]
where $T_M$ is the tangent bundle of $M$ and $N_{U/M}$ is the normal
sheaf of the local embedding $f$. In other words, $N_{U/M}$ is a local
presentation of the abelian cone stack $\NN_X$.

If $\chi:(U',M')\rightarrow(U,M)$ is a morphism of local embeddings we get an
induced commutative diagram
\[\comdia{I/I^2\resto U'}{}{{f}\upst\Omega_{M}\resto
U'}{}{}{}{I'/I'^2}{ }{{f'}\upst\Omega_{M'}\quad,}\] in other words a
homomorphism
\[\tilde{\chi}:[I/I^2\rightarrow{f}\upst\Omega_{M}]\resto
U'\longrightarrow[I'/I'^2\rightarrow{f'}\upst\Omega_{M'}]
\quad.\] We have
$\tilde{\chi}\comp\phi\resto U'=\phi'$ in $D(\O_{U'_\et})$, because of
the naturality of $\phi$. Thus the induced morphism
\[\dual{\tilde{\chi}}:[N_{U'/M'}/{f'}\upst
T_{M'}]\longrightarrow[N_{U/M}/f\upst T_M]\resto U'
 \] is compatible with the isomorphisms to
$\NN_X$. Note that, in particular,
$\dual{\tilde{\chi}}$ is an isomorphism of cone stacks over $U'$.

Recall Lemma~\ref{nci}. Let $\chi:(U',M')\rightarrow(U,M)$ be a morphism of
local embeddings. Then we get an induced morphism from the ${f'}\upst
T_{M'}$-cone $C_{U'/M'}$ to the $f\upst T_M\resto U'$-cone $C_{U/M}\resto U'$.
Note that the kernel of ${f'}\upst T_{M'}\rightarrow f\upst T_M\resto U'$ is
${f'}\upst T_{M'/M}$.

\begin{lem} The  pair $(C_{U/M}\hookrightarrow N_{U/M})\resto U'$ is the
quotient of $(C_{U'/M'}\hookrightarrow N_{U'/M'})$ by the action of ${f'}\upst
T_{M'/M}$.
\end{lem}
\begin{pf} This follows immediately from Proposition~\ref{lonc}.
\end{pf}

\begin{cor} The isomorphism
\[\dual{\tilde{\chi}}:[N_{U'/M'}/{f'}\upst
T_{M'}]\longrightarrow[N_{U/M}/f\upst T_M]\resto U'
 \] identifies the closed subcone stack $[C_{U'/M'}/{f'}\upst T_{M'}]$
with the closed subcone stack $[C_{U/M}/f\upst T_M]\resto U'$.
\end{cor}

By this corollary, there exists a unique closed subcone stack
$\CC_X\hookrightarrow \NN_X$, such that for every local embedding
$(U,M)$ of $X$ we have $\CC_X\resto U=[C_{U/M}/f\upst T_M]$, or in
other words that
\[\comdia{C_{U/M}}{}{N_{U/M}}{}{}{}{\CC_X}{}{\NN_X}\] is cartesian.

\begin{defn}
The cone stack $\CC_X$ is called the {\em intrinsic normal cone }of
$X$.
\end{defn}

\begin{them}
The intrinsic normal cone $\CC_X$ is of pure dimension zero. Its
abelian hull is $\NN_X$.
\end{them}
\begin{pf} The second claim follows because the normal sheaf is the abelian
hull of the normal cone, for any local embedding. To prove the claim
about the dimension of $\CC_X$, consider a local embedding $(U,M)$ of
$X$, giving rise to the local presentation $C_{U/M}$ of
$\CC_X$. Assume that $M$ is of pure dimension.  We then have a
cartesian and cocartesian diagram of $U$-stacks
\[\comdia{f\upst T_M\times C_{U/M}}{}{C_{U/M}}{}{}{}{C_{U/M}}{}{[C_{U/M}/f\upst
T_M].}\] Thus ${C_{U/M}}{\rightarrow}{[C_{U/M}/f\upst T_M]}$ is a
smooth epimorphism of relative dimension $\dim M$. So since $C_{U/M}$
is of pure dimension $\dim M$ (see \cite{fulton}, B.6.6) the stack
$[C_{U/M}/f\upst T_M]$ has pure dimension $\dim M-\dim M=0$.
\end{pf}

\begin{rmk}
One may construct $\NN_X$ by simply gluing the various stacks
$[N_{U/M}/f\upst T_M]$, coming from the local embeddings of $X$. So
one doesn't need the construction preceding Proposition~\ref{has} to
define the intrinsic normal sheaf and the intrinsic normal cone. But
for objects $E\com$ of $D^-(\O_{X_\et})$ satisfying ($\star$) other
than $L_X\com$, we could not prove that such gluing works a
priori. The problem is, that in general one does not have such a nice
distinguished class of local resolutions of $E\com$ (like the one
coming from local embeddings for $L_X\com$). In general, local (free)
resolutions of $E\com$ are only compatible up to homotopy.
\end{rmk}

\subsection{Basic Properties}

\begin{prop}[Local Complete Intersections] \label{colci}
The following are equivalent.
\begin{enumerate}
\item \label{colci1} $X$ is a local complete intersection,
\item \label{colci2} $\CC_X$ is a vector bundle stack,
\item \label{colci3} $\CC_X=\NN_X$.
\end{enumerate}
If, for example, $X$ is smooth, we have $\CC_X=\NN_X=BT_X$.
\end{prop}
\begin{pf}
(\ref{colci1})$\Longrightarrow$(\ref{colci3}).
If $X$ is a local complete intersection, then local embeddings of $X$
are regular immersions, but for regular immersions normal cone and
normal sheaf coincide.

(\ref{colci3})$\Longrightarrow$(\ref{colci2}).
If for a local embedding normal cone and normal sheaf coincide, then it is a
regular immersion. Thus $X$ is a local complete intersection so that
$\NN_X$ is a vector bundle stack.

(\ref{colci2})$\Longrightarrow$(\ref{colci1}).
If $\CC_X$ is a vector bundle stack it is equal to its abelian hull.
Hence $\CC_X=\NN_X$ and $X$ is a local complete intersection.
\end{pf}

\begin{prop}[Products] \label{prod}
Let $X$ and $Y$ be Deligne-Mumford stacks of finite type over
$k$. Then
\[\NN_{X\times Y}=\NN_{X}\times
\NN_{Y}\] and \[\CC_{X\times Y}=\CC_X\times \CC_Y.\]
\end{prop}
\begin{pf} If $X\subset V$ and $Y\subset W$ are affine schemes, it is easy to
check that there is a natural isomorphism $C_{X/V}\times
C_{Y/W}\to C_{X\times Y/V\times W}$, compatible with \'etale base change; the
same is true if we replace the normal cone by the normal sheaf.

If $C$ is an $E$-cone and $D$ is an $F$-cone, then $C\times D$ is an $E\times
F$-cone and there is a canonical isomorphism of cone stacks $[C/E]\times
[D/F]\to [C\times D/E\times F]$.

Putting together this remarks and verifying that the canonical isomorphisms
glue completes the proof.
\end{pf}
\begin{prop}[Pullback] \label{fshecs}
Let $f:X\to Y$ be a local complete intersection morphism. Then we have
a natural short exact sequence of cone stacks
\[\NN_{X/Y}\longrightarrow\CC_X\longrightarrow f\upst\CC_Y\]
over $X$, where $\NN_{X/Y}=h^1/h^0(T_{X/Y}\com)$.
\end{prop}
\begin{pf}
We have a distinguished triangle in $D(\O_{X_\et})$
\[\dt{f\upst L_Y}{}{L_X}{}{L_{X/Y}}{},\]
and $L_{X/Y}$ is of perfect amplitude contained in $[-1,0]$. So by
Proposition~\ref{dtscs} we have a short exact sequence of abelian cone
stacks \[\NN_{X/Y}\longrightarrow\NN_X\longrightarrow f\upst\NN_Y\] on
$X$. So the claim is local in $X$ and we may assume that we have a
diagram
\[\begin{array}{ccccc}
X & \stackrel{i}{\longrightarrow} & M'' & \longrightarrow & M' \\
 & \searrow & \downarrow & & \downarrow \\
 & & Y & \longrightarrow & M,
\end{array}\]
where the square is cartesian, the vertical maps are smooth, the
horizontal maps are local immersions, $i$ is regular and $M$ is
smooth. Then we have a morphism of short exact sequences of cones on
$X$:
\[\begin{array}{ccccc}
i\upst T_{M''/Y} & \longrightarrow & T_{M'}\resto X & \longrightarrow
& T_M\resto X \\
\downarrow & & \downarrow & & \downarrow \\
N_{X/M''} & \longrightarrow & C_{X/M'} & \longrightarrow &
C_{Y/M}\resto X.
\end{array}\]
This is a local presentation for the short exact sequence
\[\NN_{X/Y}\longrightarrow\CC_X\longrightarrow f\upst\CC_Y\]
of cone stacks.
\end{pf}

\section{Obstruction Theory}

\subsection{The Intrinsic Normal Sheaf as Obstruction}

A closed immersion $T\to \ol T$ of schemes is called a {\em
square-zero extension }with ideal sheaf $J$ if $J$ is the ideal sheaf
of $T$ in $\ol T$ and $J^2=0$.

Let $X$ be a Deligne-Mumford stack, $\NN_X$ its intrinsic normal
sheaf. Let $T\to\ol T$ be a square zero extension with ideal sheaf $J$
and $g:T\to X$ a morphism. By the functorialities of the cotangent
complex we have a canonical homomorphism
\begin{equation} \label{idkwtp}
g\upst L_X\com\longrightarrow L_T\com\longrightarrow L\com_{T/\ol T}
\end{equation}
in $D(\O_{T_\et})$.
Since $\tau_{\geq-1}L\com_{T/\ol T}=J[1]$, this homomorphism may be
considered as an element $\omega(g)$ of $\Ext^1(g\upst L_X\com,J)$.
Recall the following basic facts of deformation theory. An extension
$\ol{g}:\ol{T}\to X$ of $g$ exists if and only if $\omega(g)=0$ and if
$\omega(g)=0$ the extensions form a torsor under $\Ext^0(g\upst
L\com_X,J)=\Hom(g\upst \Omega_X,J)$.

These facts can be interpreted in terms of the intrinsic normal sheaf
$\NN_X$ of $X$. To do this, note that (\ref{idkwtp}) gives rise to a
morphism
\[h^1/h^0(L\com_{T/\ol T})\longrightarrow h^1/h^0(g\upst L\com_X)\]
of cone stacks over $T$. Since $h^1/h^0(L\com_{T/\ol T})=C(J)$ and
$h^1/h^0(g\upst L\com_X)=g\upst\NN_X$ we have constructed a morphism
$ob(g):C(J)\to g\upst\NN_X$. We also consider the morphism
$0(g):C(J)\to g\upst\NN_X$ given as the composition of $C(J)\to X$
with the vertex of $g\upst\NN_X$. By $\ul{\Hom}(ob(g),0(g))$ we shall
denote the sheaf of 2-isomorphisms of cone stacks from $ob(g)$ to
$0(g)$, restricted to $T_\et$.

Given a square zero extension $T\to\ol T$ and a morphism $g:T\to X$,
we denote the set of extensions $\ol{g}:\ol{T}\to X$ of $g$ by
$\Ext(g,\ol{T})$. These extensions in fact form a sheaf on $T_\et$
which we shall denote $\ul{\Ext}(g,\ol{T})$.

\begin{prop} \label{ince}
There is a canonical isomorphism
\[\ul{\Ext}(g,\ol{T})\longiso\ul{\Hom}_{\O_T}(ob(g),0(g))\]
of sheaves on $T_\et$. In particular, extensions of $g$ to $\ol T$
exist, if and only if $ob(g)$ is $\aaa^1$-equivariantly isomorphic to
$0(g)$.
\end{prop}
\begin{pf}
Locally, we may embed $X$ into a smooth scheme $M$ and call the
embedding $i:X\to M$, the conormal sheaf $I/I^2$. Then there always
exist local extensions $h:\ol{T}\to M$ of $i\comp g:T\to M$.
\[\comdia{T}{}{\ol T}{g}{}{h}{X}{i}{M}\]
Any such $h$ gives rise to a homomorphism $h^\sharp:g\upst I/I^2\to
J$, and hence to a realization of $ob(g)$ as the morphism of cone
stacks induced by the homomorphism of complexes
\[h^\sharp:g\upst[I/I^2\to i\upst\Omega_M]\longrightarrow[J\to0].\]
Note that if $\tilde{h}$ is another such extension, the difference
between $h$ and $\tilde h$ induces a homomorphism $g\upst i\upst
\Omega_M\to J$, which is in fact a homotopy from $h^\sharp$ to
$\tilde{h}^\sharp$.

Now let $\ol{g}:\ol{T}\to X$ be an extension of $g$. Then $(i\comp\ol
g)^\sharp=0$, so that we get a homotopy from any local $h^\sharp$ as
above to $0$, or in other words a local $\aaa^1$-equivariant
isomorphism from $ob(g)$ to $0(g)$. Since these local isomorphisms
glue, we get the required map
\[\ul{\Ext}(g,\ol{T})\longrightarrow\ul{\Hom}(ob(g),0(g)).\]
To construct the inverse, let $\theta:ob(g)\to0(g)$ be a 2-isomorphism
of cone stacks. Note that $\theta$ defines for every local $h$ as
above an extension of $h^\sharp$ to $\ol{h}^\sharp:i\upst\Omega_M\to J$
(use Lemma~\ref{lo2i}). Changing $h$ by $\ol{h}^\sharp$ defines
$h':\ol{T}\to M$ such that $(h')^\sharp=0$. Thus $h'$ factors through
$X$, and in fact these locally defined $h'$ glue to give the required
extension $\ol{g}:\ol{T}\to X$.
\end{pf}

\begin{prop}
There is a canonical isomorphism
\[\ul{\Aut}(0(g))\longiso\sheafhom(g\upst\Omega_X,J)\]
of sheaves on $T_\et$.
\end{prop}
\begin{pf}
Again, Lemma~\ref{lo2i} shows that the automorphisms of $0(g)$ are
(locally) the homomorphisms from $g\upst i\upst \Omega_M$ to $J$
vanishing on $g\upst I/I^2$. The exact sequence
\[I/I^2\longrightarrow i\upst \Omega_M\longrightarrow
\Omega_X\longrightarrow 0\]
finishes the proof. See also Lemma~\ref{lmh}.
\end{pf}

\begin{cor}
The sheaf $\ul{\Hom}(\ob(g),0(g))$ is a formal
$\sheafhom(g\upst\Omega_X,J)$-torsor. So if $ob(g)\cong0(g)$, the set
$\Hom(ob(g),0(g))$ is a torsor under the group
$\Hom(g\upst\Omega_X,J)$.
\end{cor}

\begin{note}
Combining this with Proposition~\ref{ince} gives that $\Ext(g,\ol T)$
is a $\Hom(g\upst\Omega,J)$-torsor if the obstruction vanishes,
reproving this fact from deformation theory alluded to above.
\end{note}

\subsection{Obstruction Theories}

\begin{defn} \label{doot} Let $E\com\in\ob D(\O_{X_\et})$ satisfy
($\star$) (see Definition~\ref{dost}). Then a homomorphism
$\phi:E\com\rightarrow L_X\com$ in $D(\O_{X_\et})$ is called an {\em
obstruction theory }for $X$, if $h^0(\phi)$ is an isomorphism and
$h^{-1}(\phi)$ is surjective. By abuse of language we also say that
$E\com$ is an obstruction theory for $X$.
\end{defn}

\begin{note}
By Proposition~\ref{isoiso} the homomorphism $\phi:E\com\to L\com_X$
is an obstruction theory if and only if
\[\dual{\phi}:\NN_X\longrightarrow\EE\]
is a closed immersion, where $\EE= h^1/h^0(\dual{(E_\fl\com)})$.  So
if $E\com$ is an obstruction theory and $\CC_X\subset \NN_X$ is the
intrinsic normal cone of $X$, then $\dual{\phi}(\CC_X)$ is a closed
subcone stack of $\EE$ of pure dimension zero. We sometimes call
$\dual{\phi}(\CC_X)$ the {\em obstruction cone }of the obstruction
theory $\phi:E\com\rightarrow L_X\com$.
\end{note}

Let $E\com\in\ob E(\O_{X_\et})$ satisfy ($\star$) and let
$\phi:E\com\to L\com_X$ be a homomorphism. Let
$\EE=h^1/h^0(\dual{(E_\fl\com)})$ and $\dual{\phi}:\NN_X\to\EE$ the
induced morphism of cone stacks. If $T\to\ol T$ is a square zero
extension of $k$-schemes with ideal sheaf $J$ and $g:T\to X$ is a
morphism, then we denote by $\phi\upst\omega(g)$ the image of the
obstruction $\omega(g)\in\Ext^1(g\upst L_X\com,J)$ in $\Ext^1(g\upst
E\com,J)$ and by $\dual{\phi}(ob(g))$ the composition
$$C(J)\stackrel{ob(g)}{\longrightarrow}
g\upst\NN_X \stackrel{g\upst\dual{\phi}}{\longrightarrow}g\upst\EE$$
of morphisms of cone stacks over $T$.

\begin{them}\label{ontoh1}
The following are equivalent.
\begin{enumerate}
\item \label{ont1} $\phi:E\com\to L_X\com$ is an obstruction theory.
\item \label{ont2}$\dual{\phi}:\NN_X\to\EE$ is a closed immersion of
cone stacks over $X$.
\item \label{ont3}For any $(T,\ol T, g)$ as above, the obstruction
$\phi\upst(\omega(g))\in\Ext^1(g\upst E\com,J)$ vanishes if and only
if an extension $\ol{g}$ of $g$ to $\ol T$ exists; and if
$\phi\upst(\omega(g))=0$, then the extensions form a torsor under
$\Ext^0(g\upst E\com,J)=\Hom(g\upst h^0(E\com),J)$.
\item \label{ont4}For any $(T,\ol T, g)$ as above, the sheaf of
extensions $\ul{\Ext}(g,\ol T)$ is isomorphic to the sheaf
$\ul{\Hom}(\dual{\phi}(ob(g)),0)$ of $\aaa^1$-equivariant
isomorphism from $\dual{\phi}(ob(g)):C(J)\to g\upst\EE$
to the vertex $0:C(J)\to
g\upst \EE$.
\end{enumerate}
\end{them}
\begin{pf}
The equivalence of (\ref{ont1}) and (\ref{ont2}) has already been
noted. In view of Proposition~\ref{ince} it is clear that
(\ref{ont2}) implies (\ref{ont4}). The implication
(\ref{ont4})$\Rightarrow$(\ref{ont3}) follows from Lemma~\ref{lmh}. So
let us prove that (\ref{ont3}) implies (\ref{ont1}).

To prove that $h^0(\phi)$ is an isomorphism we can assume that
$X=\spec R$ is an affine scheme (as the statement is local); let $A$
be any $R$-algebra, $M$ any $A$-module.  Let $T=\Spec A$, $\ol T=\Spec
(A\oplus M)$, where the ring structure is given by
$(a,m)(a',m')=(aa',am'+a'm)$. Let $g:T\to X$ be the morphism induced
by the $R$-algebra structure of $A$. Then $g$ extends to $\ol T$, so
there is a bijection $\Hom(h^0(L_X\com)\otimes A,M)\to
\Hom(h^0(E\com)\otimes A,M)$.
This implies easily that $h^0(\phi)$ is an isomorphism.

The fact that $h^{-1}(\phi)$ is surjective is local in the \'etale
topology (and only depends on $\tau_{\ge -1}E\com$). Assume therefore
that $X$ is an affine scheme, $i:X\to W$ a closed embedding in a
smooth affine scheme $W$, and let $I$ be the ideal of $X$ in $W$. We
can assume that $E^0=f^*\Omega_W$ (see the proof of
\ref{rmcs}), that $E^{-1}$ is a coherent sheaf, and that $E^i=0$ for $i\ne
0,-1$.

We have to prove that $E^{-1}\to I/I^2$ is surjective; let $M$ be its
image.  Let $T=X$, $\tilde M\subset I$ the inverse image of $M$, and
$\ol T\subset W$ the subscheme defined by $\tilde M$; let $g:T\to X$
be the identity. We can extend $g$ to the inclusion $\tilde g:\ol T\to
W$.  Let $\pi:I/I^2\to I/\tilde M$ be the natural projection. By
assumption $\pi$ factors via $E^0$ if and only if $g$ extends to a map
$\ol T\to X$, if and only if $\pi\circ
\phi^{-1}:E^{-1}\to I/\tilde M$ factors via $E^0$. As
$\pi\circ\phi^{-1}$ is the zero map, it certainly factors. Therefore
$\pi$ also factors. Consider now the commutative diagram with exact
rows
\[\begin{array}{ccccccc}
 E^{-1} & \longrightarrow & E^0 & \longrightarrow & h^0(E\com)
& \longrightarrow & 0 \\
\ldiag{\phi} & & \Vert  & & \Vert & & \\
 I/I^2 & \longrightarrow & E^0 & \longrightarrow & h^0(E\com) &
\longrightarrow & 0.
\end{array}\]
By an easy diagram chasing argument, the fact that $\pi$ factors via $E^0$
together with $\pi\circ\phi^{-1}=0$ implies $\pi=0$, hence
$\phi^{-1}:E^{-1}\to I/I^2$ is surjective.
\end{pf}

\noprint{
%\subsection{Intrinsic Normal Cone and Curvilinear Extensions}
A square-zero extension will be called {\em curvilinear }if it is isomorphic
to $\Spec K[t]/t^r\to \Spec K[t]/t^{r+s}$, with $K$ a field.
\begin{lem} Let $T\to \ol T$ be a curvilinear extension, $g:T\to X$ a
morphism. Then $\ob(g)$ factors via $\CC_X$. Conversely, if $\DD$ is a closed
subcone stack of $\NN_X$ such that for every curvilinear extension $T\to
\ol T$ and every morphism $g:T\to X$ the morphism $ob(g)$ factors via $\DD$,
then $\DD$ contains $\CC_X$.
\end{lem}
\begin{pf} The first statement is local, so assume that $X$ is an affine
scheme embedded as closed subscheme in a smooth scheme $W$, with ideal sheaf
$\iI$. Let $\hat T=\Spec K[[t]]$; extend $g$ to a morphism
$\hat g:\hat T\to W$. (I am not ready with it yet).
\end{pf}

}

\newcommand{\Sets}{\hbox{\sl Sets}}
\newcommand{\Art}{\hbox{\sl Art}}
\newcommand{\ann}{\hbox{\sl ann}\,}
\newcommand{\nbx}{h^1/h^0(T_X\com)}
\newcommand{\nor}{\ol{ N}_p}
\newcommand{\cp}{\ol{ C}_p}

\subsection{Obstructions for Small Extensions}

Let $\Art$ be the category of local Artinian $k$-algebras with residue
field $k$. A {\sl small extension} will be a surjective morphism $A'\to A$ in
$\Art$ with kernel $J$ isomorphic to $k$. A {\sl semi-small} extension is one
with
kernel isomorphic to a $k$-vector space as an $A'$-algebra.

Let $F:\Art \to \Sets$ be a pro-representable covariant functor (in the
sense of \cite{schlessinger}). An {\em obstruction space} for $F$ is a
set $k$-vector space $T^2$ and, for any semi-small extension $A'\to A$
with kernel $J$, an exact sequence $$ F(A')\longrightarrow
F(A)\stackrel{ob}{\longrightarrow} T^2\otimes J.$$ This means that,
for all $\xi\in F(A)$, $\xi$ is in the image of $F(A')$ if and only if
$ob(\xi)=0$. It is also required that $ob$ is functorial in the
obvious sense (see \cite{kawamata}). We say that $v\in T^2$ {\em
obstructs a small extension }$A'\to A$ if $ob(\xi)=v\otimes w$ for
some $\xi\in F(A)$ and some nonzero $w\in J$.

Let $X$ be a Deligne-Mumford stack, $p\in X$ a fixed point with
residue field $k$.  Let $h_p:\Art\to \Sets$ be the covariant functor
associating to an object $A$ of $\Art$ the set of morphisms $\spec
A\to X$ sending the closed point to $p$. The functor $h_p$ is
pro-representable, and it is unchanged if we replace $X$ by any \'etale
open neighborhood of $p$.

Let $N_p=p^*\NN_X$, and let $\nor$ be the coarse moduli space of
$N_p$. Note that $\nor =T^1_{X,p}/T^0_{X,p}$, so that $\nor$ is in
fact a $k$-vector space. Here $T^i_{X,p}=h^i(p\upst
T_X\com)=\dual{h^i(p\upst L_X\com)}$ are the `higher tangent spaces'
of $X$ at $p$. Let $\cp\subset \nor$ be the subcone coarsely
representing $p\upst\CC_X$. Proposition~\ref{ince} implies that $\nor$
is an obstruction space for $h_X$. The following is probably known but
we include a proof for lack of a suitable reference; it is a version
of Theorem~\ref{ontoh1} for semi-small extensions.

\begin{lem} The space $\nor$ is a universal
obstruction space for $h_p$; that is, for any other obstruction space $T^2$,
there is a unique injection
$N_p\to T^2$ compatible with the obstruction maps.
\end{lem}
\begin{pf} Let $(U,W)$ be a local embedding
for $X$ near
$p$. Assume that $W=\spec P$, $U=\spec R=\spec P/I$; let $\Mm$ be the maximal
ideal of $p$ in $P$, and assume that $I\subset \Mm^2$. In this case
$\nor=\dual{(I/\Mm I)}$.

If $n$ is sufficiently large, the natural map
$I/\Mm I\to (I+\Mm^n)/(\Mm I+\Mm^n)$ is an isomorphism; choose such an $n$. Let
$A'_n\to A_n$ be the extension $P/(\Mm I+\Mm^n)\to P/(I+\Mm^n)$, and let
$\xi_n\in h_p(A_n)$ be the natural quotient map. Then if $T^2$ is any
obstruction space, the obstruction to $\xi_n$ gives a linear map $I/\Mm I\to
T^2$ which must be injective. It is easy to check by functoriality that taking
a
different $n$ does not change the map. But given any semi-small extension
$A'\to
A$, there is always an extension of the type $A'_n\to A_n$ mapping to it, so
one
can apply functoriality again.
\end{pf}

\begin{prop}
Every $v\in \nor$ obstructs
some small extension; it obstructs some small curvilinear extension
if and only if
$v\in\cp$.
\end{prop}
\begin{pf}
Let $v\in \nor$, and view it as a linear map $I\to k$ having $\Mm I$ in the
kernel; we prove first that $v$ is an obstruction for some small extension. Let
$L=\ker v$, and choose $n$ sufficiently large, so that $L+\Mm^n\ne I+\Mm^n$.
Let
$A=P/I+\Mm^n$, and $A=P/L+\Mm^n$; choose $\xi:R\to A$ to be the natural
surjection. Let $J=\ker (A'\to A)$; $J$ is naturally isomorphic to $I/L$. Then
$\ob_\xi:I/\Mm I\to J$ is the obvious map, and the image of the dual map in
$\nor$ is the vector space generated by $v$.

Choose a set of generators $f_1,\ldots,f_r$ of
$I$ inducing a basis for $I/\Mm I$. This defines a map $f:W\to \aaa^r$ such
that
$U$ is the fiber over the origin. Then $\cp$ is the normal cone to the image of
$W$ in $\aaa^r$. The proof then follows the argument of Proposition 20.2 in
\cite{harris}.
\end{pf}

\section{Obstruction Theories and Fundamental Classes} \label{otafc}

\subsection{Virtual Fundamental Classes}

As usual, let $X$ be a Deligne-Mumford stack over  $k$.

\begin{defn} We call an obstruction theory $E\com\rightarrow L_X\com$ {\em
perfect}, if $E\com$ is of perfect amplitude contained in $[-1,0]$.
\end{defn}

Now assume that $X$ is separated (or, more generally, satisfies the condition
of
Vistoli in
\cite{vistoli}). We shall denote by $A_k(X)$ the rational Chow group of cycles
of dimension $k$ on
$X$ modulo rational equivalence tensored with $\qq$ (see [ibid]). We shall also
use the corresponding bivariant groups $A^k(X\rightarrow Y)$, for morphisms
$X\rightarrow Y$ of separated Deligne-Mumford stacks.

Let $E\com$ be a perfect obstruction theory for $X$, and let
$\CC_X\hookrightarrow h^1/h^0(\dual{E})$ be the intrinsic normal
cone. We call $\rk E\com$ the {\em virtual dimension }of $X$ with
respect to the obstruction theory $E\com$. Recall that $\rk E\com=\dim
E^0-\dim E^{-1}$, if locally $E\com$ is written as a complex of vector
bundles $[E^{-1}\rightarrow E^0]$. This is a well-defined locally
constant function on $X$. We shall assume that the virtual dimension
of $X$ with respect to $E\com$ is constant, equal to $n$.

To construct the {\em virtual fundamental class }$[X,E\com]\in A_n(X)$
of $X$ with respect to the obstruction theory $E\com$, we would like
to simply intersect the intrinsic normal cone $\CC_X$ with the vertex
(zero section) of $h^1/h^0(\dual{E})$. Since $h^1/h^0(\dual{E})$ is
smooth of relative dimension $-n$ over $X$, the codimension of $X$ in
$h^1/h^0(\dual{E})$ is $-n$, so that the dimension of the intersection
of $\CC_X$ with $X$ is $0-(-n)=n$. Unfortunately, this construction would
require Chow groups for Artin stacks, which we do not have at our
disposal. This is why we shall make the assumption that $E\com$ has
global resolutions.

\begin{defn} Let $F\com=[F^{-1}\rightarrow F^0]$ be a homomorphism of vector
bundles on $X$ considered as a complex of $\O_X$-modules concentrated in
degrees
$-1$ and $0$. An isomorphism $F\com\rightarrow E\com$ in $D(\O_{X_\et})$ is
called a {\em global resolution }of $E\com$.
\end{defn}

Let $F\com$ be a global resolution of $E\com$. Then
\[h^1/h^0(\dual{E})=[\dual{F^{-1}}/\dual{F^0}],\] so that $F_1=\dual{F^{-1}}$
is
a (global) presentation of $h^1/h^0(\dual{E})$. Let $C(F\com)$ be the fibered
product
\[\comdia{C(F\com)}{}{F_1}{}{}{}{\CC_X}{}{h^1/h^0(\dual{E}).}\]
Then $C(F\com)$ is a closed subcone of the vector bundle $F_1$. We
define the {\em virtual fundamental class }$[X,E\com]$ to be the
intersection of $C(F\com)$ with the zero section of $F_1$. Note that
$C(F\com)\rightarrow \CC_X$ is smooth of relative dimension $\rk F_0$
(where $F_0=\dual{F^0}$), so that $C(F\com)$ has pure dimension $\rk
F_0$ and $[X,E\com]$ then has degree
\[\rk F_0-\rk F_1=\rk E\com =n.\]

\begin{prop} \label{vfcigr}
The virtual fundamental class $[X,E\com]$ is independent of the global
resolution $F\com$ used to construct it.
\end{prop}
\begin{pf} Let $H\com$ be another global resolution of $E\com$. Without loss of
generality assume that
$H\com\rightarrow E\com$ and $F\com\rightarrow E\com$ are given by morphisms of
complexes. Then we get an induced homomorphism $H^0\oplus F^0\rightarrow E^0$.
So by constructing the cartesian diagram
\[\comdia{K^{-1}}{}{H^0\oplus F^0}{}{}{}{E^{-1}}{}{E^0,}\] and letting
$K^0=H^0\oplus F^0$, we get a global resolution $K\com$ of $E\com$ such that
both
$H\com$ and $F\com$ map to $K\com$ by a strict monomorphism. So it suffices to
compare $F\com$ with
$K\com$. Dually, we have an epimorphism $K_1\rightarrow F_1$. Consider the
diagram
\[\begin{array}{ccccc} X & \stackrel{0}{\longrightarrow} & C(H\com) &
\longrightarrow & C(F\com) \\
\ldiag{} & & \ldiag{} & & \rdiag{} \\ X & \stackrel{0}{\longrightarrow} & K_1 &
\stackrel{\alpha}{\longrightarrow} & F_1,
\end{array}\] in which both squares are cartesian. Note that $\alpha$ is
smooth.
The virtual fundamental class using $F\com$ is equal to
\[(\alpha\comp
0)\upsh[C(F\com)]=0\upsh\alpha\upsh[C(F\com)]=0\upsh[C(H\com)],\]
which is the virtual fundamental class using $H\com$.
\end{pf}

\begin{example} If $X$ is a  complete intersection, then $L_X\com$ is of
perfect amplitude contained in
$[-1,0]$, so that $L_X\com$ itself is a perfect obstruction theory.
Any embedding of $X$ into a smooth Deligne-Mumford stack gives rise to
a global resolution of $L_X\com$.The virtual fundamental class
$[X,L_X\com]$ thus obtained is equal to $[X]$, the `usual' fundamental
class.
\end{example}

\begin{numrmk}[Virtual Structure Sheaves] \label{vss}
Let $X$ be a Deligne-Mumford stack and let $\CC\hookrightarrow\EE$ be
a closed subcone stack of a vector bundle stack. Then we define a
graded commutative sheaf of coherent $\O_X$-algebras $\O_{(\CC,\EE)}$
as follows.

If $\EE\cong[E_1/E_0]$, then $\CC$ induces a cone $C$ in $E_1$ and we
set
\[\O_{(\CC,\EE)}^i=\sheaftor_i^{\O_{E_1}}(\O_C,\O_X),\]
where we think of $\O_X$ as an $\O_{E_1}$-algebra via the zero section
of $E_1$. Standard arguments show that
\[\O_{(\CC,\EE)}=\bigoplus_i \O_{(\CC,\EE)}^i\]
is independent of the choice of presentation $\EE\cong[E_1/E_0]$.
Hence the locally defined sheaves glue, giving rise to a globally
defined sheaf.

If $\CC=\CC_X$, $E\com$ is a perfect obstruction theory of $X$ and
$\EE=h^1/h^0(\dual{E\com})$, we call $\O_{(\CC,\EE)}$ the {\em virtual
structure sheaf }of $X$ with respect to the obstruction theory
$E\com$, denoted $\O_{(X,E\com)}$. This seems to be the virtual
structure sheaf proposed by Kontsevich in
\cite{K}.

If one has on $X$ a homological Chern character
$\tau:K_0(X)\rightarrow A\lst(X)$ one can define the virtual
fundamental class of $X$ with respect to $E\com$ by
\[[X,E\com]=\td(E\com)\cap\tau(\O_{(X,E\com)}).\]
This agrees with the above definition using global resolutions if they
exist. In the absence of a general Riemann Roch theorem, we rather
assume the existence of global resolutions.
\end{numrmk}

\subsection{Basic Properties}

\begin{prop}[No obstructions] \label{ehs} If $E\com$ is
perfect, $h^0(E\com)$ is locally free and $h^1({E\com})=0$, then
$X$ is smooth, the virtual dimension of $X$ with respect to $E\com$ is $\dim X$
and
the virtual fundamental class $[X,E\com]$ is just $[X]$, the usual
fundamental class. \qed
\end{prop}

\begin{prop}[Locally free obstructions] \label{xself} Let $X$ be
smooth and $E\com$ a perfect
obstruction theory for $X$. If $h^0(E\com)$ is locally free (or
equivalently $h^1(\dual{E\com})$ is locally free) then the virtual
fundamental class is
\[[X,E\com]=c_r(h^1(\dual{E\com}))\cdot[X],\]
where $r=\rk h^1(\dual{E\com})$.
\end{prop}
\begin{pf}
To see this, note that if $F\com$ is a global resolution of $E\com$,
then $C(F\com)=\im(F_0\rightarrow F_1)$.
\end{pf}

\begin{prop}[Products]
Let $E\to L_X$ be a perfect obstruction theory for $X$ and $F\to L_Y$
a perfect obstruction theory for $Y$. Then $L_{X\times Y}=L_X\boxplus
L_Y$. The induced homomorphism $E\boxplus F\to L_X\boxplus L_Y$ is a
perfect obstruction theory for $X\times Y$. If $E$ and $F$ have global
resolutions, then so does $E\boxplus F$ and we have
\[[X\times Y,E\boxplus F]=[X,E]\times[Y,F]\]
in $A_{\rk E+\rk F}(X\times Y)$.
\end{prop}
\begin{pf}
The statement about cotangent complexes is \cite{Ill}, Chapitre II, Corollaire
3.11.
To prove the rest, use Proposition~\ref{prod}.
\end{pf}

Consider a cartesian diagram of Deligne-Mumford stacks
\begin{equation}\label{asd}
\comdia{X'}{u}{X}{g}{}{f}{Y'}{v}{Y,}
\end{equation}
where $v$ is a local complete intersection morphism.
Let $E\to L_X$ and $F\to L_{X'}$ be perfect obstruction theories for
$X$ and $X'$, respectively.

\begin{defn} A {\em compatibility datum }(relative to $v$) for $E$ and
$F$ is a triple
$(\phi,\psi,\chi)$ of morphisms in $D(\O_{X'})$ giving rise to a
morphism of distinguished triangles
\[\begin{array}{ccccccc}
u\upst E & \stackrel{\phi}{\longrightarrow} & F &
\stackrel{\psi}{\longrightarrow} & g\upst L_{Y'/Y} &
\stackrel{\chi}{\longrightarrow} & u\upst E[1] \\
\ldiag{} & &  \ldiag{} & & \ldiag{} & & \rdiag{} \\
u\upst L_X & {\longrightarrow} & L_{X'} &
{\longrightarrow} & L_{X'/X} &
{\longrightarrow} & u\upst L_X[1] .
\end{array}\]
Given a compatibility datum, we call $E$ and $F$ {\em compatible
}(over $v$).
\end{defn}

Assume that $E$ and $F$ are endowed with such a compatibility datum.
Then we get (Proposition~\ref{dtscs}) a short exact sequence of vector
bundle stacks
\[g\upst h^1/h^0(T_{Y'/Y}\com)\longrightarrow
h^1/h^0(\dual{F})\longrightarrow u\upst h^1/h^0(\dual{E})\]
which we shall abbreviate by
\[g\upst\NN_{Y'/Y} \longrightarrow\FF\stackrel{\phi}{\longrightarrow}
u\upst\EE.\]

If $v$ is a regular local immersion, then $\NN_{Y'/Y}=N_{Y'/Y}$ is the
normal bundle of $Y'$ in $Y$. Its pullback to $X'$ we shall denote by
$N$.

\begin{lem} \label{ture}
If $Y$ and $Y'$ are smooth and $v$ a regular local immersion, then
there is a (canonical) rational equivalence $\beta(Y',X)\in
W\lst(N\times\FF)$ such that
\[\del\beta(Y',X)=[\phi\upst
C_{u\upst\CC_X/\CC_X}]-[N\times\CC_{X'}].\]
\end{lem}
\begin{pf}
Let $X\to M$ be a local embedding, where $M$ is smooth. We get an
induced cartesian diagram
\[\comdia{X'}{}{X}{}{}{}{Y'\times M}{}{Y\times M,}\]
which we enlarge to
\[
\begin{array}{ccccc}
N\times_XC & \longrightarrow & u\upst C & \longrightarrow & C \\
\ldiag{} & & \ldiag{} & & \rdiag{} \\
N & \longrightarrow & X' & \stackrel{u}{\longrightarrow} & X \\
\ldiag{} & & \ldiag{j} & & \rdiag{i} \\
N_{Y'/Y}\times M & \stackrel{\rho}{\longrightarrow} & Y'\times M &
\stackrel{v}{\longrightarrow} & Y\times M,
\end{array}
\]
where $C$ is the normal cone of $X$ in $Y\times M$.
As in Section~\ref{stinc} we have a canonical rational equivalence
$\beta(Y'\times M,X)\in W\lst(N\times_X C)$ such that
\[\del\beta(Y'\times M,X)=[C_{u\upst C/C}]-[N\times C_{X'/Y'\times M}].\]
By Proposition~\ref{vrqi} $\beta(Y'\times M,X)$ is invariant under the
action of $N\times u\upst i\upst T_{Y\times M} $ on $N\times_X C$. So
it descends to $N\times_X\CC_X$. In particular, $\beta(Y'\times M,X)$
is invariant under the subsheaf $N\times j\upst T_{Y'\times M}$ and
thus descends to $N\times[u\upst C/j\upst T_{Y'\times M}]$. Note that
$[u\upst C/j\upst T_{Y'\times M}]=\FF\times_{\EE}\CC_X$, which is a
closed subcone stack of $\FF$. So pushing forward via this closed
immersion, we get a rational equivalence on $N\times\FF$ which we
denote by $\beta(Y',X)$. We have
\[\del\beta(Y',X)=[\phi\upst
C_{u\upst\CC_X/\CC_X}]-[N\times\CC_{X'}]\]
as required.
Now use Proposition~\ref{vcgt} to show that $\beta(Y',X)$ does not
depend on the choice of the local embedding $X\to M$. So even if no
global embedding exists, the locally defined rational equivalences
glue, proving the lemma.
\end{pf}

\begin{prop}[Functoriality] \label{fcrp}
Let $E$ and $F$ be compatible perfect  obstruction theories, as above.
If $E$ and $F$ have global resolutions then
\[v\upsh[X,E]=[X',F]\]
holds in the following cases.
\begin{enumerate}
\item \label{fcrp1} $v$ is smooth,
\item \label{fcrp2} $Y'$ and $Y$ are smooth.
\end{enumerate}
\end{prop}
\begin{pf}
First note that one may choose global resolutions $[E_0\rightarrow
E_1]$ of $\dual{E}$ and $[F_0\to F_1]$ of $\dual{F}$ together with a
pair of epimorphisms $\phi_0:F_0\to u\upst E_0$ and $\phi_1:F_1\to
u\upst E_1$ such that
\[\comdia{F_0}{\phi_0}{u\upst E_0}{}{}{}{F_1}{\phi_1}{u\upst E_1}\]
commutes. Letting $G_i$ be the kernel of $\phi_i$ we get a short exact
sequence of homomorphisms of vector bundles
\[\begin{array}{ccccccccc}
0 & \longrightarrow & G_0 & \longrightarrow & F_0 & \longrightarrow &
u\upst E_0 & \longrightarrow & 0\\
& & \downarrow & & \downarrow & & \downarrow & & \\
0 & \longrightarrow & G_1 & \longrightarrow & F_1 & \longrightarrow &
u\upst E_1 & \longrightarrow & 0.
\end{array}\]
The induced short exact sequence
\[[G_1/G_0]\longrightarrow [F_1/F_0]\longrightarrow [u\upst E_1/u\upst
E_0]\] of vector bundle stacks is isomorphic to
$g\upst\NN_{Y'/Y}\to\FF\to \EE$. We let $C_1=\CC_X\times_{\EE}E_1$ and
$D_1=\CC_{X'}\times_{\FF}F_1$. Then $[X,E]=0_{E_1}\upsh[C_1]$ and
$[X',F]=0_{F_1}\upsh[D_1]$, where $0_{E_1}$ and $0_{F_1}$ are the zero
sections of $E_1$ and $F_1$, respectively.

If $v$ is smooth, then by Proposition~\ref{fshecs} the diagram
\[\comdia{\CC_{X'}}{}{u\upst\CC_X}{}{}{}{\FF}{}{u\upst\EE}\]
is cartesian, which implies that
\[\comdia{D_1}{}{u\upst C_1}{}{}{}{F_1}{}{u\upst E_1}\]
is cartesian. Hence $0_{u\upst E_1}\upsh[u\upst C_1]=0_{
F_1}\upsh[D_1]$ and we have
\begin{eqnarray*}
v\upsh[X,E] & = & v\upsh 0_{E_1}\upsh[C_1] \\
            & = & 0_{u\upst E_1}\upsh[u\upst C_1] \\
            & = & 0_{F_1}\upsh[D_1] \\
            & = & [X',F].
\end{eqnarray*}

If $Y'$ and $Y$ are smooth, let us first treat the case that $v$ is a
regular local immersion. Then we may choose $F_1$ as the fibered
product
\[\comdia{F_1}{}{E_1}{}{}{}{\FF}{\phi}{\EE.}\]
Lifting the rational equivalence $\beta(Y',X)$ of Lemma~\ref{ture} to
$N\times F_1$ we get that
\[[N\times D_1]=\phi\upst[C_{u\upst C_1/C_1}]\]
in $A\lst(N\times F_1)$. Then we have
\begin{eqnarray*}
[X',F_1] & = & 0_{F_1}\upsh[D_1] \\
         & = & 0_{N\times F_1}\upsh[N\times D_1] \\
         & = & 0_{N\times F_1}\upsh\phi\upst[C_{u\upst C_1/C_1}] \\
         & = & 0_{N\times u\upst E_1}\upsh[C_{u\upst C_1/C_1}] \\
         & = & 0_{\upst E_1}\upsh v\upsh[C_1] \\
         & = & v\upsh 0_{E_1}\upsh [C_1] \\
         & = & v\upsh [X,E].
\end{eqnarray*}

In the general case factor $v$ as
$$Y'\stackrel{\Gamma_v}{\longrightarrow} Y'\times Y
\stackrel{p}{\longrightarrow} Y.$$ Then Diagram~\ref{asd} factors as
\[\begin{array}{ccccc}
X' & \longrightarrow & Y'\times X & \longrightarrow & X \\
\downarrow & & \downarrow & & \downarrow \\
Y' & \stackrel{\Gamma_v}{\longrightarrow} & Y'\times Y &
\stackrel{p}{\longrightarrow} & Y.
\end{array}\]
Since $Y'$ is smooth it has a canonical obstruction theory, namely
$\Omega_{Y'}$. As obstruction theory on $Y'\times X$ take
$\Omega_{Y'}\boxplus E$. Then $\Omega_{Y'}\boxplus E$ is compatible
with $E$ over $p$ and $F$ is compatible with $\Omega_{Y'}\boxplus E$
over $\Gamma_v$. So combining Cases~(\ref{fcrp1}) and~(\ref{fcrp2})
yields the result.
\end{pf}

\section{Examples}

\subsection{The Basic Example}

Assume that
\[\comdia{X}{j}{V}{g}{}{f}{Y}{i}{W}\]
is a cartesian diagram of schemes, that $V$ and $W$ are smooth and
that $i$ is a regular embedding. Let $E\com$ be the complex
$[g\upst\dual{N_{Y/W}}\to j^*\Omega_V]$ (in degrees $-1$ and $0$),
where the map is given by pulling back to $X$ and composing
$\dual{N_{Y/W}}\to i^*\Omega_W$ with $f^*\Omega_W\to \Omega_V$. The complex
$E\com$ has a natural morphism to $L_X\com$, induced by $g^*L_Y\com\to
L_X\com$ and $j^*L_V\com\to L_X\com$ (note that $E\com$ is the
cokernel of $g^*i^*L_W\com\to j^*L_V\com\oplus g^*L_Y\com$, where the
first component is the negative of the canonical map).

This makes $E\com$ into a perfect obstruction theory for $X$; the
virtual fundamental class $[X,E\com]$ is just $i\upsh[V]$ as defined
in \cite{fulton}, p.~98. The construction also works in case $X$, $Y$,
$V$ and $W$ are assumed to be just Deligne-Mumford stacks.

\subsection{Fibers of a Morphism between Smooth Stacks}

Let $f:V\rightarrow W$ be a morphism of algebraic stacks. We shall
assume that $V$ and $W$ are smooth over $k$ and that $f$ has
unramified diagonal, so that $V$ is a relative Deligne-Mumford stack
over $W$. Let $w:\spec k\rightarrow W$ be a $k$-valued point of $W$
and let $X$ be the fiber of $f$ over $w$. In this situation $X$ has an
obstruction theory as follows.

Choose a smooth morphism $\tilde W\to W$, with $\tilde W$ a scheme,
and a lifting $\tilde w:\spec k\to \tilde W$ of $w$ (assume $k$
algebraically closed). Let $\tilde V$ be the fiber product
$V\times_W\tilde W$; by the assumptions $\tilde V$ is a smooth
Deligne-Mumford stack. Then $X$ is isomorphic to the fiber over
$\tilde w$ of $\tilde V\to \tilde W$, hence it has an obstruction
theory as above.

To check that the obstruction theory so defined does not depend on the
choices made, it is enough to compare two different ones induced by a
smooth morphism of schemes $\tilde W'\to \tilde W$; this is then a
straightforward verification. Similarly, one generalizes to the case
of arbitrary ground field $k$.

See Example~\ref{rfe} for an alternative construction.

\newcommand{\reldu}{\omega}
\newcommand{\eext}{\mathop{{\rm Ext}}\nolimits}

\subsection{Moduli Stacks of Projective Varieties}

Let $M$ and $X$ be Deligne-Mumford stacks.  Let $p:M\to X$ be a flat,
relatively
Gorenstein projective morphism: by this we mean that it has constant relative
dimension and that the relative dualizing complex $\omega\com_{M/X}$ is a line
bundle
$\omega$.

If $G\com\in D^+(\O_X)$, we have $p\upsh G\com=p\upst G\com\otimes \reldu$. So
for any complex $F\com\in D^-(\O_M)$ we have natural isomorphisms
\[\eext^k_{\O_M}(F\com,p\upst G\com)\to
 \eext^k_{\O_M}(F\com\otimes \reldu,p\upsh G\com)\to
\eext^k_{\O_X}(Rp\lst (F\com \otimes \reldu), G\com).
\]
In particular, the Kodaira-Spencer map
$L_{M/X}\to p^*L_X[1]$ induces a map
$E\com\to L\com_X$ (well-defined up to homotopy). Define the complex
$E\com$ on
$X$ to be
$Rp_*(L\com_{M/X}\otimes \reldu)[-1]$.

\begin{prop} Let $p:M\to X$ be a flat, projective, relatively Gorenstein
morphism of
Deligne-Mumford stacks, and assume that the family $M$ is universal at every
point
of $X$ (e.g., $X$ is an open set in a fine moduli space and $M$ is the
universal
family). Then $E\com\to L_X\com$ is an obstruction theory for $X$.
\end{prop}
\begin{pf}
Let $T$ be a scheme, $f:T\to X$ a morphism, and consider the cartesian
diagram \[
\comdia{N}{g}{M}{q}{}{p}{T}{f}{X.}\]
If $T\to \bar T$ is a square zero extension with ideal sheaf $\jJ$, the
obstruction to extending $N$ to a flat family over $\bar T$ lies in
$\eext^2(L\com_{N/T},q^*\jJ)$, and the extensions, if they exist, are a torsor
under $\eext^1(L\com_{N/T},q^*\jJ)$.
Now $L\com_{N/T}=g^*L\com_{M/X}$ because $p$ is flat, hence
$$\eext^k_{\O_N}(L\com_{N/T},q^*\jJ)=\eext^k_{\O_M}(L\com_{M/X},Rg_*q^*\jJ)=
\eext^k_{\O_M}(L\com_{M/X},p^*Rf_*\jJ).$$
By the
previous
argument,
$$\eext^k_{\O_M}(L\com_{M/X},p^*Rf_*\jJ)=\eext^{k-1}_{\O_X}(E\com,Rf_*\jJ)=
\eext^{k-1}_{\O_T}(f^*E\com,\jJ).$$
Assume now that $X$ is an open subset of a fine moduli space,
that is the family $M$ is universal at every point. This implies that
the fibers of $p$ have finite and reduced automorphism group, hence
$E\com$ satisfies ($\star$).

The map $E\com\to L_X\com$ induces morphisms
$$\phi_k:\eext^k_{\O_N}(L\com_{N/T},q^*\jJ)=\eext^{k-1}_{\O_T}(f^*E\com,\jJ)
\to\eext^{k-1}_{\O_T}(f^*L\com_{X},\jJ)$$ and the fact that $X$ is a moduli
space
implies that
$\phi_1$ is an isomorphism
and $\phi_2$ is
injective. By Theorem \ref{ontoh1}, this implies that $E\com$
is an obstruction theory for $X$.
\end{pf}
\begin{rmk} If $p$ is smooth of relative dimension $\le 2$, then $E\com$ is a
perfect
obstruction theory.
\end{rmk}

\subsection{Spaces of Morphisms}

Let $C$ and $V$ be projective $k$-schemes. Let $X=\Mor(C,V)$ be the $k$-scheme
of morphisms from $C$ to $V$ (see \cite{fgaIV}). Let $f:C\times X\rightarrow V$
be the universal morphism and $\pi:C\times X\rightarrow X$ the
projection. By the functorial properties of the cotangent complex we
get a homomorphism
\[f\upst L\com_V\longrightarrow L\com_{C\times X}\longrightarrow
L\com_{C\times X/C}\]
and a homomorphism
\[\pi\upst L\com_X\longrightarrow L\com_{C\times X/C}.\]
The latter is an isomorphism so that we get an induced homomorphism
\[e:f\upst L_V\com\longrightarrow\pi\upst L_X\com.\] Assume that $C$ has a
dualizing complex $\omega_C$. Then we get a homomorphism
\[e\otimes\omega_C:f\upst L_V\com\ltensor\omega_C\longrightarrow\pi\upst
L_X\com\ltensor\omega_C=\pi\upsh L_X\com\] and by adjunction a homomorphism
\[\pi\lst(e\otimes\omega_C):R\pi\lst(f\upst
L_V\com\ltensor\omega_C)\longrightarrow L_X\com.\] By duality we have
\[R\pi\lst(f\upst L_V\com\ltensor\omega_C)=\dual{(R\pi\lst(f\upst T_V\com))}.\]
Let us denote the resulting homomorphism by
\[\dual{\pi\lst(\dual{e})}:\dual{(R\pi\lst(f\upst T_V\com))}\longrightarrow
L_X\com.\]

\begin{prop}  \label{feot} Assume that $C$ is Gorenstein. Then
the homomorphism
$\phi:=\dual{\pi\lst(\dual{e})}$ is an obstruction theory for $X$. If $C$ is a
curve and $V$ is smooth then this obstruction theory is perfect.
\end{prop}
\begin{pf} Let $T$ be an affine scheme, $g:T\to X$ a
morphism, $\jJ$ a coherent sheaf on $T$; let $p:C\times T\to T$ be the
projection, $h:C\times T\to V$ the morphism induced by $g$.

By an argument analogous to that in the previous example, we get
$$\eext^k_{\O_{C\times T}}(h^*L_V\com,p^*\jJ)=
\eext^k_{\O_C}(g^*E\com,\jJ).$$
Apply now Theorem \ref{ontoh1}, more precisely the equivalence between
(1) and (3).  Choose any square zero extension $\bar T$ of $T$ with
ideal sheaf $\jJ$. Then $g$ extends to $\bar g:\bar T\to X$ if and
only if $h$ extends to $\bar h:C\times \bar T\to V$, if and only if
$\phi^*\omega(g)$ is zero in $\eext^1_{\O_{C\times
T}}(h^*L_V\com,p^*\jJ)$. The extensions, if they exist, form a torsor
under $\Hom_{\O_{C\times T}}(h^*L_V\com,p^*\jJ)$.
\end{pf}

\section{The Relative Case}

\subsection{Bivariant Theory for Artin Stacks}

For what follows, we need a little bivariant intersection theory for
algebraic stacks that are not necessarily of Deligne-Mumford type.

For simplicity, let us assume that $f:X\to Y$ is a morphism of
algebraic $k$-stacks which is representable. This assumption implies
that whenever
\[\comdia{X'}{}{Y'}{}{}{}{X}{f}{Y}\]
is a cartesian diagram and $Y'$ is a Deligne-Mumford stack satisfying
the condition needed to define its Chow group (see \cite{vistoli}),
then $X'$ is of the same type.  The following remarks can be
generalized to any morphism $f$ satisfying this property, e.g.\ any
$f$ which has finite unramified diagonal.

For such an $f:X\to Y$ we define bivariant groups $A\upst(X\to Y)$ by
using the same definition as Definition~5.1 in \cite{vistoli}. Then
just as in [ibid.] one proves that the elements of $A\upst(X\to Y)$
act on Chow groups of Deligne-Mumford stacks.

The same definition as [ibid.] Definition~3.10 applies in case $f:X\to
Y$ is a regular local immersion, and defines a canonical element
$[f]\in A\upst(X\to Y)$ whose action on cycle classes is denoted by
$f\upsh$. This is justified, since Theorems~3.11,~3.12, and~3.13 from
[ibid.] hold with the same proofs in this more general context. In
fact, $[f]$ even commutes with the Gysin morphism for any other local
regular immersion of algebraic stacks.

Similarly, if $f:X\to Y$ is flat, flat pullback of cycles defines a
canonical orientation $[f]\in A\upst(X\to Y)$.

\subsection{The Relative Intrinsic Normal Cone}

We shall now replace the base $\spec k$ by an arbitrary smooth (or
more generally pure dimensional, but always of constant dimension)
algebraic $k$-stack $Y$ (not necessarily of Deligne-Mumford type). We
shall consider algebraic stacks $X$ over $Y$ which are of relative
Deligne-Mumford type over $Y$, i.e.\ such that the diagonal
$X\rightarrow X\times_Y X$ is unramified. This assures that
$h^i(L_{X/Y}\com)=0$, for all $i>0$ (i.e.\ $h^1(L_{X/Y}\com)=0$), so
that $L_{X/Y}$ satisfies Condition~($\star$).

The {\em relative intrinsic normal sheaf }$\NN_{X/Y}$ is defined as
$$\NN_{X/Y}=h^1/h^0(T_{X/Y}\com).$$ Using local embeddings of $X$ into
schemes smooth over $Y$, we construct as in the absolute case a
subcone stack $\CC_{X/Y}\subset\NN_{X/Y}$ called the {\em relative
intrinsic normal cone }of $X$ over $Y$. If $n=\dim Y$, then
$\CC_{X/Y}$ is of pure dimension $n$.

The definition of a {\em relative obstruction theory }is the same as
Definition~\ref{doot}, with $L_X\com$ replaced by $L_{X/Y}\com$. As in
the absolute case the relative intrinsic normal cone embeds as a
closed subcone stack of a vector bundle stack
$$\CC_{X/Y}\subset h^1/h^0(\dual{E}),$$ if $E$ is a perfect relative
obstruction theory. (Note that `perfect' means `absolutely perfect'.)

So let $E$ be a perfect obstruction theory for $X$ over $Y$ admitting
global resolutions.
If $X$ is a separated Deligne-Mumford stack then we get a virtual
fundamental class $[X,E\com]\in A_{n+\rk E}(X)$ by `intersecting
$\CC_X$ with the vertex of $h^1/h^0(\dual{E})$' as in the discussion
preceding Proposition~\ref{vfcigr}.

Consider the following diagram, where $Y$ and $Y'$ are smooth of
constant dimension, $v$ has finite unramified diagonal and $X$ and
$X'$ are separated Deligne-Mumford stacks.
\begin{equation}\label{dfrcpb}
\comdia{X'}{u}{X}{}{}{}{Y'}{v}{Y}
\end{equation}

\begin{prop} \label{ncfprc}
There is a natural morphism
$$\alpha:\CC_{X'/Y'}\longrightarrow\CC_{X/Y}\times_YY'.$$ If
(\ref{dfrcpb}) is cartesian, then $\alpha$ is a closed immersion. If,
moreover, $v$ is flat, then
$\alpha$ is an isomorphism.
\end{prop}
\begin{pf}
Both statements follow immediately from the corresponding
properties of normal cones for schemes.
\end{pf}

\begin{prop}[Pullback]\label{pull}
Let $E\rightarrow L_{X/Y}$ be a perfect obstruction theory for $X$
over $Y$. If (\ref{dfrcpb}) is cartesian then $u\upst E$ is a perfect
obstruction theory for $X'$ over $Y'$. If $E$ has global resolutions
so does $u\upst E$ and for the induced virtual fundamental classes we
have
\[v\upsh[X,E]=[X',u\upst E],\]
at least in the following cases.
\begin{enumerate}
\item $v$ is flat,
\item $v$ is a regular local immersion.
\end{enumerate}
\end{prop}
\begin{pf}
Let $E^{-1}\to E^0$ be a global resolution of $E\com$ and $C$ the cone
induced by $\CC_{X/Y}$ in $E_1$. Let $u\upst E_i=E'_i$, and $D$ the
cone induced by $\CC_{X'/Y'}$ in $E'_1$.

If $v$ is flat we have $\CC_{X'/Y'}=v\upst\CC_{X/Y}$ and hence
$D=v\upst C$ by Proposition~\ref{ncfprc} and the statement follows
>from the fact that $v\upsh$ is a bivariant class; in this case that
$v\upsh$ commutes with $0\upsh_{E_1}$, where $0:X\to E_1$ is the zero
section.

If $v$ is a regular local immersion, let $N=N_{Y'/Y}$ and use
Vistoli's rational equivalence
\[\beta(Y',X)\in W\lst(N\times_Y C)\]
(see Proposition~\ref{vrqsb}) to prove that $v\upsh[C]=[D]$. Then
proceed as before.
\end{pf}

The following are relative versions of the basic properties of virtual
fundamental classes from Section~\ref{otafc}.

\begin{prop}[Locally free obstructions] \label{rulfc}
Let $E\com$ be a perfect relative obstruction theory for $X$ over $Y$
such that $h^{0}(E\com)$ is locally free. Assume that $E\com$ has
global resolutions and $X$ is a separated Deligne-Mumford stack, so
that the virtual fundamental class $[X,E\com]$ exists.
\begin{enumerate}
\item If $h^{-1}(E\com)=0$, then $X$ is smooth over $Y$ and
$[X,E\com]=[X]$.
\item If $X$ is smooth over $Y$, then $h^1(\dual{E})$ is locally free and
$[X,E\com]=c_r(h^1(\dual{E}))\cdot[X]$, where $r=\rk h^1(\dual{E})$.
\end{enumerate}
\end{prop}
\begin{pf} The proofs are the same as in the absolute case (Propositions
\ref{ehs}and \ref{xself}).
\end{pf}
\begin{prop}[Products] \label{rvopf}
Let $E$ be a perfect relative obstruction theory for $X$ over $Y$ and
$F$ a perfect relative obstruction theory for $X'$ over $Y'$. Then
$E\boxplus F$ is a perfect relative obstruction theory for $X\times
X'$ over $Y\times Y'$. If $E$ and $F$ have global resolutions and $X$
and $X'$ are separated Deligne-Mumford stacks, then $E\boxplus F$ has
global resolutions and $X\times X'$ is a separated Deligne-Mumford
stack and we have
\[[X\times X',E\boxplus F]=[X,E]\times[X',F]\]
in $A_{\dim Y+\dim Y'+\rk E+\rk F}(X\times X')$.
\end{prop}

Let $E$ be a perfect relative obstruction theory for $X$ over $Y$ and
$F$ a perfect relative obstruction theory for $X'$ over $Y$. Let
$v:Z'\rightarrow Z$ be a local complete intersection morphism of
$Y$-stacks that have finite unramified diagonal over $Y$. Let there be
given a cartesian diagram
\[\comdia{X'}{u}{X}{g}{}{f}{Z'}{v}{Z}\]
of $Y$-stacks. Then $E$ and $F$ are {\em compatible over $v$} if there
exists a homomorphism of distinguished triangles
\[\begin{array}{ccccccc}
u\upst E & {\longrightarrow} & F &
{\longrightarrow} & g\upst L_{Z'/Z} &
{\longrightarrow} & u\upst E[1] \\
\ldiag{} & &  \ldiag{} & & \ldiag{} & & \rdiag{} \\
u\upst L_{X/Y} & {\longrightarrow} & L_{X'/Y} &
{\longrightarrow} & L_{X'/X} &
{\longrightarrow} & u\upst L_{X/Y}[1] .
\end{array}\]
in $D(\O_{X'})$.

\begin{prop}[Functoriality]\label{rvofp}
If $E$ and $F$ are compatible over $v$, then
\[v\upsh[X,E]=[X',F],\]
at least if $v$ is smooth or $Z'$ and $Z$ are smooth over $Y$.
\end{prop}
\begin{pf}
The proof is the same as that of Proposition~\ref{fcrp}.
\end{pf}

\begin{numex} \label{rfe}
Consider a cartesian diagram
$$\comdia{X}{j}{V}{g}{}{h}{Y}{i}{W}$$
of algebraic stacks, where $i$ and $j$ are local immersions and $h$
has unramified diagonal. We have a canonical homomorphism
$$\phi:j\upst L_{V/W}\longrightarrow L_{X/Y},$$
which makes $j\upst L_{V/W}$ a relative obstruction theory for $X$
over $Y$. To see this, it suffices to prove that
$h^{-1}(F\com)=h^0(F\com)=0$, where $F\com$ is the cone of $\phi$. But
$F\com$ is isomorphic to the cone of the homomorphism
$$g\upst L_{Y/W}\longrightarrow L_{X/V},$$ so this is indeed true.

Now if $V$ and $W$ are smooth, then $h^i(L_{V/W})=0$ for all
$i\not=-1,0$ and $j\upst L_{V/W}$ is a perfect obstruction theory. In
particular, we get a virtual fundamental class
$$[X,j\upst L_{V/W}]\in A_{\dim Y+\dim V-\dim W}(X),$$
if $Y$ is pure dimensional, $j\upst L_{V/W}$ has global resolutions
and $X$ is a separated Deligne-Mumford stack.

If, in addition, $i$ is a regular local immersion with normal bundle
$N_{Y/W}$, the normal cone $C_{X/V}$ of $X$ in $V$ is a closed subcone
of $g\upst N_{Y/W}$ and intersecting it with the zero section $0$ of
$g\upst N_{Y/X}$ gives a class
$$0\upsh[C_{X/V}]\in A_{\dim Y+\dim V-\dim W}(X).$$
The proof that
$$0\upsh[C_{X/V}]=[X,j\upst L_{V/W}]$$
is similar to the proof of Proposition~\ref{pull}.
\end{numex}

%\bibliographystyle{/usr/local/lib/tex/inputs/plain}
%\bibliography{refer}

\begin{flushleft}
{\tt behrend@math.ubc.ca}\\
{\tt fantechi@alpha.science.unitn.it}
\end{flushleft}

\end{document}